\documentclass[10pt]{amsart}

\title[Quantum Cohomology of Weighted Projective Spaces]{The Quantum
  Orbifold Cohomology of Weighted Projective Spaces}
\author[Coates]{Tom Coates}
\address{Department of Mathematics\\
  Imperial College London\\
  Huxley Building, 180 Queen's Gate\\
  London SW7 2AZ 
  \\UK}
\email{t.coates@imperial.ac.uk}

\author[Corti]{Alessio Corti}
\address{Department of Mathematics\\
  Imperial College London\\
  Huxley Building, 180 Queen's Gate\\
  London SW7 2AZ\\
  UK}
\email{a.corti@imperial.ac.uk}

\author[Lee]{Yuan-Pin Lee}
\address{Department of Mathematics\\
  University of Utah\\
  155 South 1400 East, JWB 233\\
  Salt Lake City, Utah 84112-0090\\
  USA}
\email{yplee@math.utah.edu}

\author[Tseng]{Hsian-Hua Tseng}
\address{Department of Mathematics\\
  University of British Columbia\\
  1984 Mathematics Road\\
  Vancouver, B.C. V6T 1Z2\\
  Canada}
\curraddr{Department of Mathematics\\
  University of Wisconsin--Madison\\
  Van Vleck Hall\\ 
  480 Lincoln Drive\\
  Madison, WI 53706-1388\\
  USA}
\email{tseng@math.wisc.edu}

\usepackage{hyperref}

\usepackage[initials,short-publishers]{amsrefs}
\usepackage{amsfonts}
\usepackage{amssymb}
\usepackage[all]{xy}

\BibSpec{article}{%
+{}{\sc \PrintAuthors} {author}
+{,}{ \textrm} {title}
+{,}{ \textit} {journal}
+{,}{ } {volume}
+{}{ \parenthesize} {date}
+{,}{ } {pages}
+{.}{} {transition}
+{}{ } {status}
+{.}{} {transition}
+{}{ \tt} {eprint}
+{.}{} {transition}
}

\BibSpec{book}{%
+{}{\sc \PrintAuthors} {author}
+{,}{ \textit} {title}
+{.}{ } {series}
+{,}{ } {volume}
+{.}{ } {publisher}
+{,}{ } {place}
+{,}{ } {date}
+{.}{} {transition}
}

\BibSpec{collection.article}{
+{} {\sc \PrintAuthors} {author}
+{,} { \textrm } {title}
+{,} { in \it \PrintConference} {conference}
+{,} { pp.~} {pages}
+{.} {\PrintBook} {book}
+{,} { \PrintDateB} {date}
+{.} {} {transition}
}

\BibSpec{innerbook}{
+{.} { \emph} {title}
+{.} { } {part}
+{:} { \emph} {subtitle}
+{.} { } {series}
+{,} { } {volume}
+{.} { Edited by \PrintNameList} {editor}
+{.} { Translated by \PrintNameList}{translator}
+{.} { \PrintContributions} {contribution}
+{.} { } {publisher}
+{.} { } {organization}
+{,} { } {address}
+{,} { \PrintEdition} {edition}
+{,} { \PrintDateB} {date}
+{.} { } {note}
}

\newcommand{\stackX}{\mathbb{P}^\mathbf{w}}
\newcommand{\PP}{\mathbb{P}}
\newcommand{\wP}{\PP(w_0,\ldots,w_n)}
\newcommand{\CC}{\mathbb{C}}
\newcommand{\ZZ}{\mathbb{Z}}
\newcommand{\RR}{\mathbb{R}}
\newcommand{\QQ}{\mathbb{Q}}
\newcommand{\TT}{\mathbb{T}}
\newcommand{\Cstar}{\CC^\times}
\newcommand{\cD}{\mathcal{D}}
\newcommand{\cI}{\mathcal{I}}
\newcommand{\inertia}{\mathcal{I}\stackX}
\newcommand{\Horb}{H_{\text{\rm orb}}^\bullet(\stackX;\CC)}
\newcommand{\HorbL}{H_{\text{\rm orb}}^\bullet(\stackX;\Lambda)}
\newcommand{\HorbLinv}{H_{\text{\rm orb}}^\bullet(\stackX;\Lambda[Q^{-1}])}
\newcommand{\cO}{\mathcal{O}}
\newcommand{\ev}{\mathrm{ev}}
\newcommand{\LX}{L\stackX}

\newcommand{\tLX}{\widetilde{\LX}}
\newcommand{\HF}{HF^\bullet(\LX)}
\newcommand{\eqHF}{HF^\bullet_{S^1}(\LX)}

\newcommand{\iso}{\cong}
\newcommand{\NovZ}{\CC(\!(Q)\!)}
\newcommand{\NovzZ}{\CC[z](\!(Q)\!)}
\newcommand{\bP}{\mathbf{P}}
\newcommand{\bQ}{\mathbf{Q}}
\newcommand{\bJ}{\mathbf{J}}
\newcommand{\Lpol}{L_{\text{\rm poly}}}
\newcommand{\Mbar}{\overline{\mathcal{M}}}
\renewcommand{\(}{\left(}
\renewcommand{\)}{\right)}
\DeclareMathOperator{\Loc}{Loc}
\newcommand{\fun}{\mathbf{1}}
\DeclareMathOperator{\age}{age}
\DeclareMathOperator{\orbdeg}{orbdeg}
\newcommand{\cC}{\mathcal{C}}
\newcommand{\cM}{\mathcal{M}}
\newcommand{\clM}{\overline{\mathcal{M}}}
\newcommand{\cU}{\mathcal{U}}
\newcommand{\cX}{\mathcal{X}}
\newcommand{\cIX}{\mathcal{IX}}
\newcommand{\HorbX}{H^\bullet_{\text{\rm orb}}(\cX;\CC)}
\newcommand{\HorbXL}{H^\bullet_{\text{\rm orb}}(\cX;\Lambda)}

\newcommand{\correlator}[1]{\left \langle #1 \right \rangle}
\newlength{\mybracketspacing}
\newcommand{\Correlator}[1]{
  \settowidth{\mybracketspacing}{$\big\langle$} 
  \big \langle \hspace{-1.3\mybracketspacing} \big \langle 
  #1 
  \big \rangle \hspace{-1.3\mybracketspacing} \big \rangle^{\cX}_\tau
}
\newcommand{\BigCorrelator}[1]{
  \settowidth{\mybracketspacing}{$\left\langle\right.$} 
  \left \langle \hspace{-1.2\mybracketspacing} \left \langle 
      #1 
    \right \rangle \hspace{-1.2\mybracketspacing} \right \rangle^{\cX}_\tau}
\DeclareMathOperator{\id}{id}

\theoremstyle{plain}

\newtheorem{thm}{Theorem}[section]
\newtheorem{pro}[thm]{Proposition}
\newtheorem{lem}[thm]{Lemma}
\newtheorem{cla}[thm]{Claim}
\newtheorem{cor}[thm]{Corollary}

\theoremstyle{definition}

\newtheorem{dfn}[thm]{Definition}
\newtheorem{nt}[thm]{Notation}

\newtheorem{rem}[thm]{Remark}

\DeclareMathOperator{\Mor}{Mor}
\DeclareMathOperator{\Aut}{Aut}

\DeclareMathOperator{\Fix}{Fix}

\DeclareMathOperator{\hcf}{hcf}

\DeclareMathOperator{\Hom}{Hom}

\DeclareMathOperator{\Spec}{Spec}

\DeclareMathOperator{\Ext}{Ext}

\renewcommand{\o}{\mathcal{O}}

\renewcommand{\r}{\mathcal{R}}

\newcommand{\p}{\mathbb{P}}
\newcommand{\C}{\mathbb{C}}
\newcommand{\z}{\mathbb{Z}}

\begin{document}

\begin{abstract}
  We calculate the small quantum orbifold cohomology of arbitrary
  weighted projective spaces. We generalize Givental's heuristic
  argument, which relates small quantum cohomology to
  $S^1$-equivariant Floer cohomology of loop space, to weighted
  projective spaces and use this to conjecture an explicit formula for
  the small $J$-function, a generating function for certain genus-zero
  Gromov--Witten invariants. We prove this conjecture using a method
  due to Bertram. This provides the first non-trivial example of a
  family of orbifolds of arbitrary dimension for which the small
  quantum orbifold cohomology is known.  In addition we obtain
  formulas for the small $J$-functions of weighted projective complete
  intersections satisfying a combinatorial condition; this condition
  naturally singles out the class of orbifolds with terminal
  singularities.
\end{abstract}

\maketitle

\tableofcontents

\section{Introduction}
\label{sec:intro}

In this paper we calculate the small quantum orbifold cohomology ring
of weighted projective space $\stackX = \wP$.  Our approach is
essentially due to Givental \cite{Givental:homological, Givental:ICM,
  Givental:equivariant}.  We begin with a heuristic argument relating
the quantum cohomology of $\stackX$ to the $S^1$-equivariant Floer
cohomology of the loop space $L\stackX$, and from this conjecture a
formula for a certain generating function --- the \emph{small
  $J$-function} --- for genus-zero Gromov--Witten invariants of
$\stackX$.  The small $J$-function determines the small quantum
orbifold cohomology of $\stackX$.  We then prove that our conjectural
formula for the small $J$-function is correct by analyzing the
relationship between two compactifications of the space of
parametrized rational curves in $\stackX$: a toric compactification
(which is closely related to our heuristic model for the Floer
cohomology of $L \stackX$) and the space of genus-zero stable maps to
$\stackX \times \PP(1,r)$ of degree ${1 \over r}$ with respect to the
second factor.  These compactifications carry natural
$\Cstar$-actions, which one can think of as arising from rotation of
loops, and there is a map between them which is $\Cstar$-equivariant.
Our formula for the small $J$-function can be expressed in terms of
integrals of $\Cstar$-equivariant cohomology classes on the toric
compactification.  Following Bertram \cite{Bertram}, we use
localization in equivariant cohomology to transform these into
integrals of classes on the stable map compactification.  This
establishes our formula for the small $J$-function, and so allows us
to determine the small quantum orbifold cohomology ring of $\stackX$.

We now give precise statements of our main results.  The reader
unfamiliar with orbifolds or with quantum orbifold cohomology may wish
first to read Section~\ref{sec:orb_coh}, where various basic features
of the theory are outlined.  Let $w_0,\ldots,w_n$ be a sequence of
positive integers and let $\stackX$ be the weighted projective space
$\wP$, \emph{i.e.} the quotient
\[
\Big[ \left(\CC^{n+1} - \{0\}\right) / \Cstar \Big]
\]
where $\Cstar$ acts with weights $-w_0,\ldots,-w_n$.  Components of
the inertia stack of $\stackX$ correspond to elements of the set
\[
F = \left\{{\textstyle {k \over w_i}} \mid 0 \leq k < w_i, \; 0 \leq i
\leq n\right\}
\]
via
\[
\inertia = \coprod_{f \in F} \, \PP(V^f),\] where $\PP(V^f)$ is the
locus of points of $\stackX$ with isotropy group containing $\exp(2
\pi \sqrt{-1} f) \in \Cstar$.  This locus is itself a weighted projective
space, of dimension
\[
\dim_f = \# \left \{j \colon w_j f \in \ZZ \right \} - 1.
\]
The orbifold cohomology $\Horb$ is equal as a vector space to
\[
\bigoplus_{f \in F} H^\bullet(\PP(V^f);\CC).
\]
It carries two ring structures and two gradings: the usual cup product
on the cohomology of $\inertia$, the Chen--Ruan orbifold cup product,
the usual grading on the cohomology of $\inertia$, and a grading where
the degree of a cohomology class is shifted by a rational number (the
{\em degree-shifting number} or {\em age}) depending on the component
of $\inertia$ on which the class is supported. In this paper, unless
otherwise stated, all products should be taken with respect to the
orbifold cup product; the degree of an element of $\Horb$ always
refers to its age-shifted degree. The involution $\zeta \mapsto
\zeta^{-1}$ on $\Cstar$ induces an involution $I$ on $\inertia$ which
exchanges $\p(V^f)$ with $\p(V^{1-f})$, $f\not = 0$, and is the
identity on $\p(V^0)$.

Since $\PP(V^0)=\stackX$, there is a canonical inclusion
$H^\bullet(\stackX;\CC) \subset \Horb$.  Let $P \in
H^2_{\text{orb}}(\stackX;\CC)$ be the image of $c_1(\cO(1))\in
H^2(\stackX;\CC)$ under this inclusion and let $Q$ be the generator
for $H_2(\stackX;\CC)$ dual to $c_1(\cO(1))$.  For each $f \in F$,
write $\fun_f$ for the image of $\fun \in H^\bullet(\PP(V^f);\CC)$
under the inclusion $H^\bullet(\PP(V^f);\CC) \subset \Horb$.  We will
often work with orbifold cohomology with coefficients in the ring
\[
\Lambda = \CC[\![Q^{1/\mathrm{lcm}(w_0,\ldots,w_n)}]\!].
\] 
This plays the role of the Novikov ring (see \cite{Manin}*{III 5.2.1}
and \cite{Hofer--Salamon}) in the quantum cohomology of
manifolds\footnote{If we were being more careful, we could take the
  Novikov ring to be the semigroup ring $R$ of the semigroup of
  degrees of effective possibly-stacky curves in $\stackX$.  But the
  degree of such a curve is $k/\mathrm{lcm}(w_0,\ldots,w_n)$ for some
  integer $k$, and so $R$ is naturally a subring of $\Lambda$.}.  The
quantum orbifold cohomology of $\stackX$ is a family of
$\Lambda$-algebra structures on $\HorbL$ parameterized by $\Horb$.
When the parameter is restricted to lie in $H^2(\stackX;\CC) \subset
\Horb$, we refer to the resulting family of algebras as the
\emph{small quantum orbifold cohomology} of $\stackX$.

Let $f_1,\ldots,f_k$ be the elements of $F$ arranged in increasing
order, and set $f_{k+1} = 1$.  The classes
\begin{equation} \label{1eq:obviousbasis}
  \begin{split}
    &\fun_{f_1}, \fun_{f_1} P, \ldots, \fun_{f_1} P^{\dim_{f_1}}, \\
    &\fun_{f_2}, \fun_{f_2} P, \ldots, \fun_{f_2} P^{\dim_{f_2}}, \\
    &\ldots, \\
    &\fun_{f_k}, \fun_{f_k} P, \ldots, \fun_{f_k} P^{\dim_{f_k}}
  \end{split}
\end{equation}
form a $\Lambda$-basis for $\HorbL$.

\begin{thm}[see Corollary~\ref{5cor:theorem1.1}]
\label{1thm:1}
  The matrix, with respect to the above basis, of multiplication by
  the class $P$ in the small quantum orbifold cohomology algebra of
  $\stackX = \wP$ corresponding to the point $tP \in H^2(\stackX;\CC)$
  is
\[
\begin{pmatrix}
0 & 0 & 0 & \cdots & 0 & r_N \\
r_1 & 0 & 0 & \cdots & 0 & 0 \\
0 & r_2 & 0 & \ddots & & 0 \\
\vdots & & \ddots & & \vdots & \vdots \\
0 & & & \ddots & 0 & 0 \\
0 & 0 & \cdots & 0 & r_{N-1} & 0 \\
\end{pmatrix}
\]
where 
\begin{align*}
  N &= \dim_{f_1} + \ldots + \dim_{f_k} + k; \\
  r_i &=
  \begin{cases}
    Q^{f_{j+1} - f_j} \displaystyle {s_{j+1} \over s_j} 
    {e^{f_{j+1} t} \over e^{f_j t}} &
    \parbox{3in}{if $i = \dim_{f_1} + \ldots + \dim_{f_j} + j$\\
      \phantom{\qquad} for some $j \leq k$;}
    \\
    1 & \text{otherwise} \\
  \end{cases} \\
  \intertext{and}
  s_j &= 
  \begin{cases}
    1 & j=1 \\
    \prod_{i=0}^n w_i^{-\left \lceil f_j w_i \right \rceil}
    & 2 \leq j \leq k+1.
  \end{cases}
\end{align*}
\end{thm}

The underlined superscript here denotes a falling factorial:
\[
x^{\underline{n}} = x (x-1) (x-2) \cdots (x-n+1).
\]

\begin{cor}
  \label{1cor:1}
  The small quantum orbifold cohomology algebra of $\stackX$ is
  the free $\Lambda$-module which is generated as a $\Lambda$-algebra
  by the classes
  \[
  \fun_{f_1}, \fun_{f_2},\ldots,\fun_{f_k} \quad \text{and} \quad P
  \] 
  with identity element $\fun_{f_1} = \fun_0$ and relations generated by
  \begin{align}
    \label{1eq:present}
    P^{\dim_{f_j} + 1} \fun_{f_j} &=
    Q^{f_{j+1} - f_j} {e^{f_{j+1}t} s_{j+1} \over e^{f_j t} s_j} \, \fun_{f_{j+1}} ,
    && 1 \leq j \leq k.
  \end{align}
  In particular, 
  \[
  P^N = {Q e^t \over w_0^{w_0} w_1^{w_1} \cdots w_n^{w_n}} \,
  \fun_0.
  \]
  If we invert $Q$ then the small quantum orbifold cohomology algebra
  is generated by $P$.  
\end{cor}

\begin{rem}
  If we set $Q$ to zero in \eqref{1eq:present} then we obtain a
  presentation for the Chen--Ruan orbifold cohomology ring of
  $\stackX$.
\end{rem}

\begin{rem}
  The combinatorial factors $r_i$ and $s_j$ can be simplified by
  rescaling the basis \eqref{1eq:obviousbasis}, replacing $\fun_f$ by
  $s_f e^{f t} \fun_f$.  See Section~\ref{sec:qu_coh_wps} for a
  precise statement.
\end{rem}

\begin{rem} \label{1rem:HF}
  Multiplication by $P$ preserves the $\CC[\![Q]\!]$-submodule of
  $\HorbL$ with basis 
  \begin{equation} \label{1eq:basis}
    \begin{split}
      &Q^{f_1} \fun_{f_1}, Q^{f_1} \fun_{f_1} P, \ldots, Q^{f_1} \fun_{f_1}
      P^{\dim_{f_1}}, \\
      &Q^{f_2} \fun_{f_2}, Q^{f_2} \fun_{f_2} P, \ldots, Q^{f_2} \fun_{f_2}
      P^{\dim_{f_2}}, \\
      &\ldots, \\ 
      &Q^{f_k} \fun_{f_k}, Q^{f_k} \fun_{f_k} P, \ldots, Q^{f_k} \fun_{f_k}
      P^{\dim_{f_k}}.
    \end{split}
  \end{equation}
  We will see in Section~\ref{sec:floer} that, after inverting $Q$, we
  can think of this submodule as the Floer cohomology of the loop
  space $\LX$.
\end{rem}

\begin{rem}
  Theorem~\ref{1thm:1} and Corollary~\ref{1cor:1} confirm the
  conjectures of Etienne Mann \cite{Mann}.  In the case of
  $\PP(w_0,w_1)$, we recover the result of \cite[Section~9]{AGV:2}.
  The Chen--Ruan orbifold cohomology ring of weighted projective
  space, which is obtained from the quantum cohomology ring by setting
  $Q=0$, has been studied by a number of authors.  Weighted projective
  space is a toric Deligne--Mumford stack --- this is spelled out in
  \cite{Boissiere--Mann--Perroni} --- so one can compute the orbifold
  cohomology ring using results of Borisov--Chen--Smith
  \cite{Borisov--Chen--Smith}.  One can also apply the methods of
  Chen--Hu \cite{Chen--Hu}, Goldin--Holm--Knutson
  \cite{Goldin--Holm--Knutson}, or Jiang \cite{Jiang}.  The
  relationship between the orbifold cohomology ring of certain
  weighted projective spaces and the cohomology ring of their crepant
  resolutions has been studied by Boissiere--Mann--Perroni
  \cite{Boissiere--Mann--Perroni}.  The relationship between the
  \emph{quantum} orbifold cohomology ring of certain weighted
  projective spaces and that of their crepant resolutions is
  investigated in \cite{CCIT:wallcrossing}.
\end{rem}

The small $J$-function of $\stackX$, a function of $t\in \CC$ taking
values in $\HorbL\otimes\CC(\!(z^{-1})\!)$, is defined by:
\[
J_{\stackX}(t) = z \, e^{P t/z}  \( 1 + \sum_{\substack{d \colon d > 0  \\
    \langle d \rangle \in F}} Q^d e^{dt} \,
(I \circ \ev_1)_\star \left( \fun_{0,1,d}^{\text{vir}}\cap { 1 \over z(z -
    \psi_1)} \right) \).
\]
Here $\fun^\text{vir}_{0,1,d}$ is the virtual fundamental class of the
moduli space $\stackX_{0,1,d}$ of genus-zero one-pointed stable maps
to $\stackX$ of degree $d$; the degree of a stable map is the integral
of the pull-back of the K\"ahler class $P$ over the domain curve;
$\langle d \rangle = d - \lfloor d \rfloor$ denotes the fractional
part of the rational number $d$; $\ev_1\colon\stackX_{0,1,d} \to
\inertia$ is the evaluation map at the marked point\footnote{This
  evaluation map does not in fact exist, but one can to all intents
  and purposes pretend that it does. See the discussion in
  Section~\ref{2sec:eval}.}; $\psi_1$ is the first Chern class of the
universal cotangent line at the marked point; and we expand the
expression $(z - \psi_{1})^{-1}$ as a power series in $z^{-1}$.  Note
that the degrees $d$ occurring in the sum will in general be
non-integral.  We will see in Section~\ref{sec:orb_coh} below that the
small $J$-function determines the small quantum orbifold cohomology of
$\stackX$: it satisfies a system of differential equations whose
coefficients are the structure constants of the small quantum orbifold
cohomology algebra.

\begin{thm}[see Corollary~\ref{4cor:1}] \label{1thm:smallJ}
  The small $J$-function $J_{\stackX}(t)$ is equal to
  \[
  z \, e^{P t/z} \sum_{\substack{d \colon d\geq 0\\ \langle d\rangle\in F}}
  {Q^d e^{dt} \over
    \prod_{i=0}^n\prod_{\substack{b \colon \langle b \rangle = \langle dw_i \rangle \\
        0<b\leq dw_i}}(w_iP+bz)} \fun_{\langle d\rangle}.
  \]
\end{thm}

From this, we deduce
\begin{cor}
  The small $J$-function $J_{\stackX}(t)$ satisfies the differential
  equation
\[
\prod_{i=0}^{n} \, \prod_{k=0}^{w_i-1} \left( w_i z {\partial
\over
\partial t} - k z\right ) \, J_{\stackX}(t) = Q e^t \, J_{\stackX}(t).
\]
\end{cor}

\subsection*{Weighted Projective Complete Intersections}
\label{sec:weight-proj-compl}

Let $\cX$ be a quasismooth complete intersection in $\p^\mathbf{w}$ of
type $(d_0,d_1, \dots, d_m)$ and let $\iota \colon \cX \to \stackX$ be
the inclusion. Define
\begin{equation}
  \label{1eq:I}
  I_{\cX}(t)=
  z \, e^{P t/z} \sum_{\substack{d \colon d\geq 0\\ \langle d\rangle\in F}}
  Q^d e^{dt}
  { \prod_{j=0}^m\prod_{\substack{b \colon \langle b \rangle = \langle
        d d_j \rangle \\
        0\leq b\leq dd_i}}(d_j P+bz)
    \over
    \prod_{i=0}^n\prod_{\substack{b \colon \langle b \rangle = \langle dw_i \rangle \\
        0<b\leq dw_i}}(w_iP+bz)} \fun_{\langle d\rangle}.
\end{equation}

\begin{cor}[see Section~\ref{sec:weighted_intersections}]
  \label{1cor:ci_specific}
  Let $k_{\cX} = \sum_{j=0}^m d_j -\sum_{i=0}^n w_i$ and let
  \begin{align}
    \label{1eq:sharp_box_condition}
    k_f &=\sum_{j=0}^m\lceil fd_j\rceil -\sum_{i=0}^n \lceil
    fw_i\rceil \\&= k_{\cX} f +\sum_{j=0}^m\langle -fd_j\rangle
    -\sum_{i=0}^n \langle -fw_i\rangle. \notag
  \end{align}
  Suppose that for each non-zero $f\in F$ we have either
  $k_f < -1$ or 
  \[
  \# \left\{ j \mid d_j f \in \ZZ \right\}
  \geq 
  \# \left\{ i \mid w_i f \in \ZZ \right\}.
  \]
  Then:
  \begin{enumerate}
  \item if $k_{\cX} <-1$ then 
\[
I_{\cX}(t)= \iota_\star \(z+tP+O(z^{-1})\)
\]
and
\[
\iota_\star \, J_{\cX} (t)=I_{\cX}(t);
\]
\item if $k_{\cX} =-1$ then 
\[
I_{\cX}(t)= \iota_\star \(z+tP+s(t)\fun_0+O(z^{-1})\)
\]
where $s(t)= Qe^t (\prod_{j=0}^m d_j !)/(\prod_{i=0}^n w_i !)$, and 
\[
\iota_\star \(e^{{s(t)}/{z}} J_{\cX} (t)\)=I_{\cX}(t);
\]
\item if $k_{\cX}=0$ then 
\[
I_{\cX}(t)= \iota_\star \(F(t)z+g(t)P+O(z^{-1})\)
\]
for some functions $F\colon \C \to \Lambda$, $g \colon \C \to \Lambda$, and 
\[
\iota_\star \, J_{\cX} \(\tau(t)\)=\frac{I_{\cX}(t)}{F(t)}
\]
where the change of variables $\tau(t)=g(t)/F(t)$ is invertible.
  \end{enumerate}
\end{cor}

The assumptions of Corollary~\ref{1cor:ci_specific} have a geometric
interpretation:

\begin{pro}[see Section~\ref{sec:weighted_intersections}]
  \label{1pro:terminal_singularities}
  The following conditions on $\cX$ are equivalent:
  \begin{enumerate}
  \item $\cX$ is well-formed and has terminal singularities.
  \item For all non-zero $f \in F$, either $
    \# \left\{ j \mid d_j f \in \ZZ \right\}
    \geq 
    \# \left\{ i \mid w_i f \in \ZZ \right\}
    $ or 
    \begin{equation}
      \label{1eq:condition_box}
      \sum_{i=0}^n \langle fw_i\rangle > 1+\sum_{j=0}^m \langle fd_j\rangle.
    \end{equation}
\end{enumerate}

In particular, if $k_{\cX} \leq 0$ and $\cX$ has terminal
singularities then the assumptions of Corollary~\ref{1cor:ci_specific}
are satisfied. If $\cX$ is Calabi--Yau then these assumptions are
equivalent to $\cX$ having terminal singularities.
\end{pro}

\begin{rem}
  We were surprised to discover the notion of terminal singularities
  occurring so naturally in Gromov--Witten theory.
\end{rem}

\begin{rem}
  \label{1rem:ci_conditions} Corollary~\ref{1cor:ci_specific}
  determines the part of the small $J$-function of $\cX$ involving
  classes pulled back from $\stackX$, and hence the part of the small
  quantum orbifold cohomology algebra of $\cX$ generated by such
  classes. 
\end{rem}

\begin{rem}
  Corollary~\ref{1cor:ci_specific} is an immediate consequence of a
  more general result, Corollary~\ref{6cor:hyp} below, which is
  applicable to any weighted projective complete intersection $\cX$
  with $k_\cX \geq 0$ and which determines the part of the ``big
  $J$-function'' of $\cX$ involving classes pulled back from
  $\stackX$.  The big $J$-function is defined in
  Section~\ref{3sec:J_FUN}.
\end{rem}

\begin{rem}
  \label{1rem:calabi_yau}
  In dimension 3, a Calabi--Yau orbifold has terminal singularities if
  and only if it is smooth. Thus Corollary~\ref{1cor:ci_specific}
  applies to only 4 of the 7555 quasismooth Calabi--Yau 3-fold weighted
  projective hypersufaces\footnote{See Gavin Brown's graded ring
    database \tt{http://www.kent.ac.uk/ims/grdb/}}:
  \begin{align*}
    X_5 &\subset \p(1,1,1,1,1)\\
    X_6 &\subset \p(1,1,1,1,2)\\
    X_8 &\subset \p(1,1,1,1,4)\\
    X_{10} &\subset \p(1,1,1,2,5).
  \end{align*}
  These can all be handled using methods of Givental
  \cite{Givental:equivariant} and others, by resolving the
  singularities of the ambient space.  In dimension $4$, however,
  there are many Gorenstein terminal quotient singularities and
  consequently many interesting examples.  For instance,
  \[
  X_7 \subset \PP(1,1,1,1,1,2)
  \]
  can be treated using Corollary~\ref{1cor:ci_specific} but not, to
  our knowledge, by existing methods.
\end{rem}

\begin{rem}
  \label{rem:1}
  Let $\cX \subset \p^\mathbf{w}$ be a quasismooth hypersurface of
  degree $d=\sum_{i=0}^nw_i$. The $I$-function of $\cX$ is a fundamental
  solution of the ordinary differential equation:
  \begin{equation}
    \label{1eq:mirrors}
H^{red} I = 0 
\quad
\text{where}
\quad
H=\prod_{i=0}^n \prod_{k=0}^{w_i-1} \Bigl(w_i
\frac{\partial}{\partial t} -k \Bigr)
-Qe^t \prod_{k=0}^{d-1} \Bigl(d\frac{\partial}{\partial t} -k \Bigr)    
  \end{equation}
  and the superscript ``red'' means that we are taking the main
  irreducible constituent: the operator obtained by removing factors
  that are common to both summands. It is shown in
  \cite[Theorem~1.1]{Corti--Golyshev} that the local system of
  solutions of Equation~(\ref{1eq:mirrors}) is
  $\text{gr}^W_{n-1}R^{n-1}f_!  \mathbb{R}_Y$, where $f\colon Y \to
  \Cstar$ is the \emph{mirror-dual Landau--Ginzburg model}:
\[
Y=
\begin{cases}
  \prod_{i=0}^ny_i^{w_i}=t\\
  \sum_{i=0}^ny_i = 1
\end{cases}
\subset
\(\Cstar\)^{n+1} \times \Cstar.
\]
This is a mirror theorem for quasismooth Calabi--Yau weighted
projective hypersurfaces.
\end{rem}

\subsection*{Acknowledgements}
\label{sec:acknowledgements}

We would like to thank Hiroshi Iritani for many stimulating
conversations, and for pointing out a number of errors in earlier
versions of this paper.  We thank the referee for many useful
comments.  A.C. had several helpful conversations about aspects of
this project with Martin Guest.  In particular, Martin suggested to
use Birkhoff factorization to recover quantum cohomology from the
$J$-function; this works, but we preferred to adopt a more explicit
approach here.  The project owes a great deal to Alexander Givental
who, directly or indirectly, taught us much of what we know about this
subject.

T.C. was partially supported by the Clay Mathematics Institute, the
Royal Society, NSF grant DMS-0401275, and a postdoctoral fellowship at
the Mathematical Sciences Research Institute.  Y.-P.L. was partially
supported by the NSF and by an AMS Centennial Fellowship.  H.-H.T. was
partially supported by a postdoctoral fellowship at the Mathematical
Sciences Research Institute.

\section{Orbifold Cohomology and Quantum Orbifold Cohomology}
\label{sec:orb_coh}

In this section we give an introduction to the cohomology and quantum
cohomology of orbifolds following \cite{AGV:1, AGV:2}.  An alternative
exposition can be found in \cite{Tseng}.  We work in the algebraic
category, using the term ``orbifold'' to mean ``smooth separated
Deligne--Mumford stack of finite type over $\CC$''.  Gromov--Witten
theory for orbifolds was originally constructed in the symplectic
setting by Chen and Ruan \cite{Chen--Ruan:orbifold, Chen--Ruan:GW}.
Note that we do not require our orbifolds to be reduced (in the sense
of Chen and Ruan): the stabilizer of the generic point of an orbifold
is allowed to be non-trivial.

\subsection{Orbifold Cohomology}

Let $\cX$ be a stack. Its inertia stack $\cIX$ is the fiber product
\[
\xymatrix{
\cIX \ar[r] \ar[d] & \cX \ar[d]^\Delta \\
\cX\ar[r]^\Delta & \cX \times \cX \\ 
}
\]
where $\Delta$ is the diagonal map.  The fiber product is taken in the
$2$-category of stacks.  One can think of a point of $\cIX$ as
a pair $(x,g)$ where $x$ is a point of $\cX$ and $g \in
\Aut_{\cX}(x)$.  There is an involution $I\colon\cIX \to
\cIX$ which sends the point $(x,g)$ to $(x,g^{-1})$.

The {\em Chen--Ruan orbifold cohomology groups} $\HorbX$ of a
Deligne--Mumford stack $\cX$ are the cohomology groups of its inertia
stack\footnote{An introduction to the cohomology of stacks can be
  found in Section~2 of \cite{AGV:2}.}:
\[
\HorbX = H^\bullet(\cIX; \CC).
\]
If $\cX$ is compact then there is an inner product, the {\em
  orbifold Poincar\'e pairing}, on orbifold cohomology defined by
\begin{align*}
  \HorbX \otimes \HorbX
  & \longrightarrow \CC \\
  \alpha \otimes \beta & \longmapsto \int_{\cIX} \alpha \cup I^\star
  \beta.
\end{align*}
We denote the pairing of $\alpha$ and $\beta$ by
$(\alpha,\beta)_{\text{orb}}$.

To each component $\cX_i$ of the inertia stack $\cIX$
we associate a rational number, the {\em age} of $\cX_i$,
defined as follows.  Choose a geometric point $(x,g)$ of
$\cX_i$ and write the order of $g \in \Aut_{\cX}(x)$
as $r$. The automorphism $g$ acts on the tangent space $T_x
\cX$, so we can write
\[
T_x \cX = \bigoplus_{0 \leq j < r} E_j
\]
where $E_j$ is the subspace of $T_x \cX$ on which $g$ acts by
multiplication by $\exp(2\pi\sqrt{-1}j/r)$.  The age of
$\cX_i$ is
\[
\age\(\cX_i\) = \sum_{j=0}^{r-1} {j \over r} \dim E_j.
\]
This is independent of the choice of geometric point $(x,g) \in
\cX_i$.  

We use these rational numbers to equip the orbifold cohomology
$\HorbX$ with a new grading: if $\alpha \in
H^p(\cX_i;\CC) \subset \HorbX$ then
the \emph{orbifold degree} or \emph{age-shifted degree} of $\alpha$ is
\[
\orbdeg(\alpha) = p + 2 \age(\cX_i).
\]
Note that $(\alpha,\beta)_{\text{orb}} \ne 0$ only if $\orbdeg \alpha +
\orbdeg \beta = 2 \dim_\CC \cX$, so for a compact orbifold
$\cX$ the orbifold cohomology
$\HorbX$ is a graded inner product space.

Weighted projective space $\stackX$ is the stack quotient 
\begin{equation}
  \label{2eq:defofwps}
  \Big[ \( {\CC^{n+1} - \{0\}} \) / \Cstar \Big]
\end{equation}
where $\Cstar$ acts with weights $-w_0,\ldots,-w_n$.  As discussed in
Section~\ref{sec:intro}, components of the inertia stack of $\stackX$
are indexed by
\[
F = \left\{{\textstyle {k \over w_i}} \, \Big\vert \, \, 0 \leq k < w_i, \; 0 \leq i
\leq n\right\}
\]
via
\[
\inertia = \coprod_{f \in F} \, \PP(V^f);
\] 
here 
\[
V^f = \left\{\(x_0,\ldots,x_n\) \in \CC^{n+1} \mid \text{$x_i = 0$
  unless $w_i f \in \ZZ$} \right\}
\]
and $\PP(V^f) = \left[ \( {V^f - \{0\}} \) / \Cstar \right]$, so that
$\PP(V^f)$ is the locus of points of $\stackX$ with isotropy group
containing $\exp(2 \pi \sqrt{-1} f) \in \Cstar$.  The involution $I$
maps the component $\PP(V^f)$ to the component $\PP(V^{\langle -f
  \rangle})$.  The age of $\PP(V^f) \subset \inertia$ is $\langle -
w_0 f \rangle + \cdots + \langle - w_n f \rangle$ where, as before,
$\langle r \rangle$ denotes the fractional part $r - \lfloor r
\rfloor$ of $r$.

\begin{rem}
  It is logical to take the action of $\Cstar$ on $\CC^{n+1}$ to have
  \emph{negative} weights $-w_0$, \ldots, $-w_n$, as we now explain.
  One could repeat all discussions in this paper working equivariantly
  with respect to the (ineffective) action of the torus $\TT^{n+1}$ on
  $\stackX$.  This action descends from  an action of $\TT^{n+1}$ on
  $\CC^{n+1}$, and we should choose this action so that
  $H^0(\stackX,\cO(1))$ is the standard representation of
  $\TT^{n+1}$.  This means that $\TT^{n+1}$ acts with \emph{negative}
  weights on $\CC^{n+1}$:
  \[
  (t_0,\ldots,t_n): (x_0,\ldots,x_n) \mapsto (t_0^{-1}
  x_0,\ldots,t_n^{-1} x_n).
  \]
  The action of $\Cstar$ in \eqref{2eq:defofwps} is obtained from the
  $\TT^{n+1}$-action on $\CC^{n+1}$ via the map
  \begin{align*}
    \Cstar \to \TT^{n+1}, && t \mapsto (t^{w_0},\ldots,t^{w_n}), 
  \end{align*}
  and so the weights of the $\Cstar$-action on $\CC^{n+1}$ should be
  negative.  To obtain the results which hold if the $\Cstar$-action
  in \eqref{2eq:defofwps} is taken with \emph{positive} weights
  $w_0,\ldots,w_n$, the reader should just replace the class $\fun_f$
  with the class $\fun_{\langle -f \rangle}$ throughout
  Section~\ref{sec:intro}.
\end{rem}

\begin{rem}
  One could instead define the orbifold cohomology of a
  Deligne--Mumford stack $\cX$ to be the cohomology of its
  \emph{cyclotomic inertia stack} constructed in
  \cite[Section~3.1]{AGV:2}, or as the cohomology of its
  \emph{rigidified cyclotomic inertia stack}
  \cite[Section~3.4]{AGV:2}.  Geometric points of the cyclotomic
  inertia stack are given by representable morphisms $B\mu_r \to \cX$.
  The rigidified cyclotomic inertia stack is obtained from the
  cyclotomic inertia stack by removing the canonical copy of $\mu_r$
  from the automorphism group of each component parametrizing
  morphisms $B\mu_r \to \cX$: this process is called
  ``rigidification'' \cite{Abramovich--Corti--Vistoli}.  From the
  point of view of calculation, it does not matter which definition
  one uses.  With our definitions,
  \[
  \PP(V^f) = \PP(w_{i_1},\ldots,w_{i_m})
  \]
  where $w_{i_1},\ldots,w_{i_m}$ are the weights $w_j$ such that $w_j
  f \in \ZZ$.  The reader who prefers the cyclotomic inertia stack ---
  which has the advantage that its components are parameterized by
  representations, and one can define the age of a representation
  without choosing a preferred root of unity --- should take
  \[
  \PP(V^f) = \PP(w_{i_1},\ldots,w_{i_m})
  \]
  but regard the index $f$ not as the rational number $\frac{j}{r}$
  (in lowest terms) but as the character $\zeta \mapsto \zeta^j$ of
  $\mu_r$.  The reader who prefers the rigidified cyclotomic
  inertia stack should similarly regard $f$ as a character of
  $\mu_r$, but take
  \[
  \PP(V^f) = \PP\({w_{i_1} \over r},\ldots,{w_{i_m} \over r}\).
  \]
\end{rem}

\subsection{Ring Structures on Orbifold Cohomology}

The orbifold cup product and the quantum orbifold product are defined
in terms of Gromov--Witten invariants of $\cX$.  These invariants are
intersection numbers in stacks of twisted stable maps to $\cX$.

\subsubsection{Moduli Stacks of Twisted Stable Maps}

Recall \cite[Section~4]{AGV:2} that an \emph{$n$-pointed twisted
  curve} is a connected one-dimensional Deligne--Mumford stack such
that:
\begin{itemize}
\item its coarse moduli space is an $n$-pointed pre-stable curve: a
  possibly-nodal curve with $n$ distinct smooth marked points;
\item it is a scheme away from marked points and nodes;
\item it has cyclic quotient stack structures at marked points;
\item it has \emph{balanced} cyclic quotient stack structures at nodes: near
  a node, the stack is \'etale-locally isomorphic to
  \[
  \Big[ \big( \Spec \CC[x,y]/(xy) \big) / \mu_r \Big]
  \]
  where $\zeta \in \mu_r$ acts as $\zeta\colon(x,y) \mapsto (\zeta x,
  \zeta^{-1} y)$.
\end{itemize}

A family of $n$-pointed twisted curves over a scheme $S$ is a flat
morphism $\pi\colon\cC \to S$ together with a collection of $n$ gerbes over
$S$ with disjoint embeddings into $\cC$ such that the geometric fibers
of $\pi$ are $n$-pointed twisted curves.  Note that the gerbes over
$S$ defined by the marked points need not be trivial: this will be
important when we discuss evaluation maps below.

An $n$-pointed \emph{twisted stable map} to $\cX$ of genus $g$ and degree $d
\in H_2(\cX;\QQ)$ is a representable morphism $\cC \to \cX$ such that:
\begin{itemize}
\item $\cC$ is an $n$-pointed twisted curve;
\item the coarse moduli space $C$ of $\cC$ has genus $g$;
\item the induced map of coarse moduli spaces $C \to X$ is stable in
  the sense of Kontsevich \cite{Kontsevich};
\item the push-forward $f_\star [\cC]$ of the fundamental class of
  $\cC$ is $d$.
\end{itemize}

A family of such objects over a scheme $S$ is a family of twisted
curves $\pi\colon\cC \to S$ together with a representable morphism $\cC \to
\cX$ such that the geometric fibers of $\pi$ give $n$-pointed twisted
stable maps to $\cX$ of genus $g$ and degree $d$.  The moduli stack
parameterizing such families is called the \emph{stack of twisted
  stable maps to $\cX$}.  It is a proper Deligne--Mumford stack, which
we denote by $\cX_{g,n,d}$.  In \cite{AGV:1,AGV:2} a very similar
object is denoted by $\mathcal{K}_{g,n}(\cX,\beta)$: the only
difference is that Abramovich--Graber--Vistoli take the degree $\beta$
to be a curve class on the coarse moduli space of $\cX$ whereas we
take $d$ to lie in $H_2(\cX;\QQ)$.  When we specialize to the case of
weighted projective space we will identify degrees $d \in
H_2(\stackX;\QQ)$ with their images under the isomorphism $
H_2(\stackX;\QQ) \iso \QQ$ given by cap product with $c_1(\cO(1))$.

\subsubsection{Evaluation Maps} \label{2sec:eval}

Given an $n$-pointed twisted stable map $f\colon\cC \to \cX$, each marked
point $x_i$ determines a geometric point $(f(x_i),g)$ of the inertia
stack $\cIX$ where $g$ is defined as follows.  Near $x_i$, $\cC$ is
isomorphic to $[\CC/\mu_r]$ and since $f$ is representable it
determines an injective homomorphism $\mu_r \to \Aut_{\cX}
\(f(x_i)\)$.  We work over $\CC$, so we have a preferred
generator $\exp(2 \pi \sqrt{-1}/r)$ for $\mu_r$.  The automorphism $g$ is
the image of this generator in $\Aut_{\cX}\(f(x_i)\)$.  Thus each
marked point gives an evaluation map to $\cIX$ defined on geometric
points of $\cX_{g,n,d}$.

These maps do \emph{not} in general assemble to give maps of stacks
$\cX_{g,n,d} \to \cIX$.  This is because things can go wrong in
families: given a family
\[
\xymatrix{
\cC \ar[r]^f \ar[d]^{\pi} & \cX \\
S & \\ }
\]
of twisted stable maps, each marked point determines a $\mu_r$-gerbe
over $S$ (for some $r$) and this gerbe will map to the inertia stack
only if it is trivial.  But, as is explained carefully in
\cite{AGV:2}, there \emph{are} evaluation maps to the rigidified
cyclotomic inertia stack and one can use this to define push-forwards
\[
\(\ev_i\)_\star \colon H^\bullet(\cX_{g,n,d};\CC) \longrightarrow \HorbX
\]
and pull-backs
\[
\(\ev_i\)^\star \colon \HorbX \longrightarrow H^\bullet(\cX_{g,n,d};\CC) 
\]
which behave as if evaluation maps $\ev_i\colon\cX_{g,n,d} \to \cIX$
existed.  We will write as if the maps $\ev_i$ themselves existed,
referring to ``the image of $\ev_i$'' etc.  This is an abuse of
language, but no ambiguity should result.

\subsubsection{Gromov--Witten Invariants}

The stack $\cX_{g,n,d}$ can be equipped \cite[Section~4.5]{AGV:2} with
a virtual fundamental class in $H_\bullet(\cX_{g,n,d};\CC)$.  In
general, $\cX_{g,n,d}$ is disconnected and its virtual dimension ---
the homological degree of the virtual fundamental class --- is
different on different components.  On the substack
$\cX_{g,n,d}^{i_1,\ldots,i_n}$ of twisted stable maps such that the
image of $\ev_m$ lands in the component $\cX_{i_m}$ of the inertia
stack, the \emph{real} virtual dimension is
\begin{equation}
  \label{2eq:dim}
2 n + (2-2g)(\dim_{\CC} \cX - 3) - 2 K_{\cX}(d) - 2 \sum_{m=1}^{n}
\age\(\cX_{i_m}\).   
\end{equation}
We will write $\(\stackX\)_{g,n,d}^{f_1,\ldots,f_n}$ for the substack of
$\stackX_{g,n,d}$ consisting of twisted stable maps such that the $m$th
marked point maps to the component $\PP(V^{f_m})$ of $\inertia$, and
denote the virtual fundamental class of $\stackX_{g,n,d}$ by $\fun^{vir}_{g,n,d}$.

There are line bundles
\begin{align*}
  L_i \to \cX_{g,n,d} && i \in \{1,2,\ldots,n\},
\end{align*}
called \emph{universal cotangent lines}, such that the fiber of $L_i$
at the stable map $f\colon\cC \to \cX$ is the cotangent line to the
coarse moduli space of $\cC$ at the $i$th marked point.  We denote the
first Chern class of $L_i$ by $\psi_i$.  There is a canonical map from
$\cX_{g,n,d}$ to the moduli stack $X_{g,n,d}$ of stable maps to the
coarse moduli space $X$ of $\cX$; the bundle $L_i$ is the pull-back to
$\cX_{g,n,d}$ of the $i$th universal cotangent line bundle on
$X_{g,n,d}$.

\emph{Gromov--Witten invariants} are intersection numbers of the form
\begin{equation}
  \label{2eq:GW}
  \int_{\cX_{g,n,d}^{vir}} \prod_{i=1}^n \ev_i^\star \alpha_i \cdot
  \psi_i^{k_i} 
\end{equation}
where $\alpha_1,\ldots,\alpha_n \in \HorbX$;
$k_1,\ldots,k_n$ are non-negative integers; and the integral means cap
product with the virtual fundamental class.  If any of the $k_i$ are
non-zero then \eqref{2eq:GW} is called a \emph{gravitational
  descendant}.  We will use correlator notation for Gromov--Witten
invariants, writing \eqref{2eq:GW} as
\[
\correlator{\alpha_1 \psi_1^{k_1},\ldots,\alpha_n
  \psi_n^{k_n}}^{\cX}_{g,n,d}.
\]

\begin{rem}
  One could avoid the complications caused by the non-existence of the
  maps $\ev_i$ by defining orbifold cohomology in terms of the
  rigidified cyclotomic inertia stack: evaluation maps to this flavour
  of inertia stack certainly exist.  Or one could replace
  $\cX_{g,n,d}$ with a moduli stack of stable maps with sections to
  all gerbes.  We will do neither of these things.  In each case there
  is a price to pay: to get the correct Gromov--Witten invariants ---
  the invariants which participate in the definition of an associative
  quantum product --- one must rescale all virtual fundamental classes
  by rational numbers depending on the stack structures at marked
  points.  This is described in detail in \cite[Section~1.4]{AGV:2}
  and \cite{Tseng}.
\end{rem}

\subsubsection{The Orbifold Cohomology Ring}

The Chen--Ruan orbifold cup product
$\underset{\scriptscriptstyle CR}{\cup}$ is defined by
\[
\( \alpha \underset{\scriptscriptstyle CR}{\cup} \beta, \gamma
\)_{\text{orb}} = \correlator{\alpha,\beta,\gamma}^{\cX}_{0,3,0}
\]
It gives a super-commutative and associative ring structure on
orbifold cohomology, called the \emph{orbifold cohomology ring}.  As
indicated in Section~\ref{sec:intro}, unless otherwise stated all
products of orbifold cohomology classes are taken using this ring
structure.

\subsubsection{Quantum Orbifold Cohomology}

Quantum orbifold cohomology is a family of $\Lambda$-algebra
structures on $\HorbXL$, where $\Lambda$ is an
appropriate Novikov ring, defined by
\begin{equation}
  \label{2eq:qcp}
  \(\alpha \bullet_\tau \beta,\gamma\)_{\text{orb}} = 
  \sum_d \sum_{n \geq 0} {Q^d \over n!}
  \correlator{\alpha,\beta,\gamma,\tau,\tau,\ldots,\tau}^{\cX}_{0,n+3,d}.
\end{equation}
Here the first sum is over degrees $d$ of effective possibly-stacky
curves in $\cX$, and $Q^d$ is the element of the Novikov ring
corresponding to the degree $d \in H_2(\cX;\QQ)$.  In the case $\cX =
\stackX$, where $H_2(\cX;\QQ)$ is one-dimensional and
\[
\Lambda = \CC[[Q^{1/\mathrm{lcm}(w_0,\ldots,w_n)}]],
\]
the element of $\Lambda$ corresponding to $d \in H_2(\cX;\QQ)$ is
$Q^{\int_d c_1(\cO(1))}$.  To interpret \eqref{2eq:qcp}, choose a
basis $\phi_1,\ldots,\phi_N$ for $\HorbX$ and set 
\[
\tau = \tau^1 \phi_1 + \cdots + \tau^N \phi_N.
\]
Then the right-hand side of \eqref{2eq:qcp} is a formal power series
in $\tau^1,\ldots, \tau^N$ and so \eqref{2eq:qcp} defines a family of
product structures $\bullet_\tau$ parameterized by a formal
neighbourhood of zero in $\HorbX$.  The WDVV
equations \cite{AGV:2, Chen--Ruan:GW} imply that this is a family of
associative products.

\emph{Small quantum orbifold cohomology} is the family $\circ_\tau$ of
$\Lambda$-algebra structures on $\HorbXL$ defined by restricting the
parameter $\tau$ in $\bullet_\tau$ to lie in a formal neighbourhood of
zero in $H^2(\cX;\CC) \subset \HorbX$.  The family is entirely
determined by its element at $\tau=0$.  This follows from the Divisor
Equation \cite{AGV:2}*{Theorem~8.3.1}:
\[
\correlator{\alpha_1,\ldots,\alpha_n,\gamma}^{\cX}_{0,n+1,d} = 
\(\int_d \gamma\)
\correlator{\alpha_1,\ldots,\alpha_n}^{\cX}_{0,n,d} 
\]
whenever $\gamma \in H^2(\cX;\CC)$ and either $d \ne 0$ or $n \geq 3$.
For example in the case $\cX = \stackX$, if $P$ is the first Chern
class of $\cO(1)$ and $t$ lies in a formal neighbourhood of zero in
$\CC$ then
\begin{equation}
  \label{2eq:smallQC}
  \( \alpha \circ_{t P} \beta, \gamma \)_{\text{orb}} = 
  \sum_{d \geq 0} Q^d e^{d t}
  \correlator{\alpha,\beta,\gamma}^{\stackX}_{0,3,d}.
\end{equation}
Analogous statements hold for general $\cX$.

\subsection{The $J$-Function}
\label{3sec:J_FUN} 
Let us write 
\[
\Correlator{\alpha_1 \psi_1^{i_1},\ldots,\alpha_m
  \psi_m^{i_m}} = \sum_d \sum_{n \geq 0} {Q^d \over n!}
\correlator{\alpha_1 \psi_1^{i_1},\ldots,\alpha_m
  \psi_m^{i_m},\tau,\tau,\ldots,\tau}^{\cX}_{0,m+n,d},
\]
so that
\[
\(\alpha \bullet_\tau \beta, \gamma\)_{\text{orb}} =
\Correlator{\alpha,\beta,\gamma}.
\]
The \emph{$J$-function of $\cX$} is
\begin{equation}
  \label{2eq:J}
  \bJ_{\cX}(\tau) = z + \tau + \BigCorrelator{\phi^\epsilon \over z - \psi_1} \phi_\epsilon,
\end{equation}
where $\phi^1,\ldots,\phi^N$ is the basis for $\HorbX$ such that
$\(\phi^i,\phi_j\)_{\text{orb}} = \delta^i_{\phantom{i}j}$; here and
henceforth we use the summation convention, summing over repeated
indices, and expand $(z - \psi_1)^{-1}$ as a power series in $z^{-1}$.
The $J$-function is a function of $\tau \in \HorbX$ taking values in
$\HorbXL \otimes \CC(\!(z^{-1})\!)$, defined for $\tau$ in a formal
neighbourhood of zero.  In other words, just as for \eqref{2eq:qcp},
we regard the right-hand side of \eqref{2eq:J} as a formal power
series in the co-ordinates $\tau^1,\ldots,\tau^N$ of $\tau$.  To
distinguish it from the small $J$-function of $\cX$ defined below,
we will sometimes refer to $\bJ_\cX$ as the \emph{big $J$-function} of
  $\cX$.

\begin{lem}
  The $J$-function satisfies 
  \begin{equation}
    \label{2eq:diffeq}
    z {\partial \over \partial \tau^i} {\partial \over \partial
      \tau^j} \bJ_{\cX} (\tau)
     = c(\tau)_{ij}^{\phantom{ij}\mu} {\partial \over \partial
       \tau^\mu} \bJ_{\cX} (\tau)
  \end{equation}
  where
  \[
  \phi_i \bullet_\tau \phi_j = c(\tau)_{ij}^{\phantom{ij}\mu} \phi_\mu.
  \]
\end{lem}

\begin{proof}
  This follows from the \emph{topological recursion relations}
  \begin{align*}
    \Correlator{\alpha \psi_1^{a+1}, \beta \psi_2^b, \gamma
      \psi_3^c} = 
    \Correlator{\alpha \psi_1^a,\phi_\mu} \,
    \Correlator{\phi^\mu, \beta \psi_2^b, \gamma
      \psi_3^c}, 
    && a,b,c \geq 0,
  \end{align*}
  exactly as in \cite{Pandharipande}.  A proof of the topological
  recursion relations is sketched in \cite{Tseng}.  For
  \begin{align*}
    z {\partial \over \partial \tau^i} {\partial \over \partial
      \tau^j} \bJ_{\cX} (\tau) 
    &= \sum_{m \geq 0} {1 \over z^m}
    \Correlator{\phi^\epsilon \psi_1^m,\phi_i,\phi_j}
    \phi_\epsilon \\
    &= \Correlator{\phi^\epsilon,\phi_i,\phi_j} \phi_\epsilon
     + \sum_{m \geq 1} {1 \over z^m} 
     \Correlator{\phi^\epsilon \psi_1^{m-1},\phi_\mu}
     \Correlator{\phi^\mu,\phi_i,\phi_j}
     \phi_\epsilon \\
    &= \Correlator{\phi_i,\phi_j,\phi^\mu} {\partial \over \partial
      \tau^\mu} \bJ_{\cX} (\tau) \\
    \intertext{and}
    \phi_i \bullet_\tau \phi_j &= \Correlator{\phi_i,\phi_j,\phi^\mu}
    \phi_\mu.
  \end{align*}
\end{proof}

The $J$-function determines the quantum orbifold product, as
\begin{equation}
  \label{2eq:Jtoqc}
  z {\partial \over \partial \tau^i} {\partial \over \partial
    \tau^j} \bJ_{\cX} (\tau)
  =
  \phi_i \bullet_\tau \phi_j + O(z^{-1}).
\end{equation}

\subsubsection{The Small $J$-Function} \label{2sec:smallJ} The small
$J$-function $J_{\cX}(\tau)$ is obtained from the $J$-function
$\bJ_{\cX}(\tau)$ by restricting $\tau$ to lie in a formal neighbourhood of
zero in $H^2(\cX;\CC) \subset \HorbX$.  In the case
of weighted projective space, we regard the small $J$-function as
being defined on a formal neighbourhood of zero in $\CC$, setting
\[
J_{\stackX}(t) = \bJ_{\stackX}(t P).
\]
\begin{lem}
  \[
  J_{\stackX}(t) = z \, e^{P t/z} \(1 + \sum_{d > 0} Q^d e^{d t}
  \correlator{ \phi^\epsilon \over z(z - \psi_1)}^{\stackX}_{0,1,d} \phi_\epsilon\).
  \]
\end{lem}

\begin{proof}
  This follows from the Divisor Equation \cite{AGV:2}*{Theorem~8.3.1}:
  \begin{multline*}
    \correlator{\alpha_1 \psi_1^{i_1},\ldots,\alpha_n
      \psi_n^{i_n},\gamma}^{\cX}_{0,n+1,d} = \(\int_d \gamma\)
    \correlator{\alpha_1 \psi_1^{i_1},\ldots,\alpha_n
      \psi_n^{i_n}}^{\cX}_{0,n,d} \\
    + \sum_{j=1}^n \correlator{\alpha_1 \psi_1^{i_1},\ldots, (\alpha_j
      \gamma) \psi_j^{i_j-1},\ldots, \alpha_n
      \psi_n^{i_n}}^{\cX}_{0,n,d}
  \end{multline*}
  whenever $\gamma \in H^2(\cX;\CC)$ and either $d \ne 0$ or $n \geq 3$.  We have
  \begin{equation}
    \label{2eq:Jexplicit}
    J_{\stackX}(t) = z + t P + 
    \sum_{\substack{d : d \geq 0 \\  \langle d \rangle \in F}} 
    \sum_{n \geq 0} \sum_{m \geq 0} {Q^d t^n \over n! \, z^{m+1}} 
    \correlator{\phi^\epsilon \psi_1^m,P,P,\ldots,P}^{\stackX}_{0,n+1,d} \phi_\epsilon.
  \end{equation}
  Now, using the Divisor Equation,
  \begin{multline} \label{2eq:smallJmainterms}
    \sum_{\substack{d : d > 0 \\ \langle d \rangle \in F}} \sum_{n
      \geq 0} \sum_{m \geq 0} {Q^d t^n \over n! \, z^{m+1}}
    \correlator{\phi^\epsilon
      \psi_1^m,P,P,\ldots,P}^{\stackX}_{0,n+1,d}
    \phi_\epsilon \\
    \begin{aligned}
      &= \sum_{\substack{d : d > 0 \\ \langle d \rangle \in F}}
      \sum_{n \geq 0} \sum_{m \geq 0} {Q^d t^n \over n! \, z^{m+1}}
      \correlator{\phi^\epsilon \psi_1^m \big( \textstyle {P \over z}
        + d \big)^n }^{\stackX}_{0,1,d}
      \phi_\epsilon \\
      & = \sum_{\substack{d : d > 0 \\ \langle d \rangle \in F}} Q^d
      \correlator{e^{P t/z} e^{d t} \phi^\epsilon \over z -
        \psi_1}^{\stackX}_{0,1,d}
      \phi_\epsilon  \\
      & = z \, e^{P t/z} \sum_{\substack{d : d > 0 \\ \langle d
          \rangle \in F}} Q^d e^{d t} \correlator{\phi^\epsilon \over
        z(z - \psi_1)}^{\stackX}_{0,1,d} \phi_\epsilon.
    \end{aligned}
  \end{multline}
  The terms in \eqref{2eq:Jexplicit} which are not in
  \eqref{2eq:smallJmainterms} are
  \[
  z + t P + 
  \sum_{n \geq 2} \sum_{m \geq 0} {t^n \over n! \, z^{m+1}} 
  \correlator{\phi^\epsilon \psi_1^m,P,P,\ldots,P}^{\stackX}_{0,n+1,0}
  \phi_\epsilon.
  \]
  Using the Divisor Equation again, this is
  \begin{equation}
    \label{2eq:remaining}
    z + t P + 
    \sum_{n \geq 2} \sum_{m \geq 0} {t^n \over n! \, z^{m+1}} 
    \correlator{\phi^\epsilon P^{n-2} \psi_1^{m-n+2},P,P}^{\stackX}_{0,3,0}
    \phi_\epsilon
  \end{equation}
  and since $L_1$ is trivial on $(\stackX)_{0,3,0}$ the summand
  vanishes unless $m = n-2$.  Thus \eqref{2eq:remaining} is
  \begin{align*}
    z + t P + 
    \sum_{n \geq 2} {t^n \over n! \, z^{n-1}} 
    \(\phi^\epsilon P^{n-2} \underset{\scriptscriptstyle CR}{\cup}
    P,P\)_{\text{\rm orb}}
    \phi_\epsilon 
    & = z + t P + 
    \sum_{n \geq 2} {t^n P^n \over n! \, z^{n-1}} \\
    & = z \, e^{P t/z}.
  \end{align*}
  Combining this with \eqref{2eq:smallJmainterms} gives
  \[
  J_{\stackX}(t) = z \, e^{P t/z} \(1 + \sum_{d > 0} Q^d e^{d t}
  \correlator{ \phi^\epsilon \over z(z - \psi_1)}^{\stackX}_{0,1,d} \phi_\epsilon\).
  \]
\end{proof}

From \eqref{2eq:Jtoqc}, we see that the small quantum cohomology
algebra is determined by
\begin{align*}
  \left. {\partial \bJ_{\cX} \over \partial \tau^j}(\tau)
  \right|_{\tau \in H^2(\cX;\CC) \subset \HorbX} && j \in \{1,2,\ldots,N\}.
\end{align*}
Let $v,w \in \HorbX$ and let $\nabla_v$ denote the directional
derivative along $v$, so that 
\[
\nabla_v \bJ(\tau) = v^\alpha {\partial \bJ \over \partial
  \tau^\alpha}(\tau)
\]
where $v = v^1 \phi_1 + \cdots + v^N \phi_N$ and $\tau = \tau^1 \phi_1
+ \cdots + \tau^N \phi_N$.  From \eqref{2eq:diffeq},
\begin{align*}
  z \, \nabla_v \nabla_w \, \bJ_\cX(\tau) &= \nabla_{v \bullet_\tau w}
  \, \bJ(\tau) \\
  &= v \bullet_\tau w + O(z^{-1}).
\end{align*}
Taking $\tau \in H^2(\cX;\CC) \subset \HorbX$ gives
\begin{equation}
  \label{2eq:smallQDEs}
  \begin{split}
    z \, \nabla_v \nabla_w \, \bJ_\cX(\tau) &= \nabla_{v \circ_\tau w}
    \, \bJ(\tau) \\
    &= v \circ_\tau w + O(z^{-1}),
  \end{split}
\end{equation}
and it follows that the small $J$-function determines the subalgebra
of the small quantum orbifold cohomology algebra which is generated by
$H^2(\cX;\CC)$.  We will see below that for weighted projective spaces
this subalgebra is the whole of the small quantum orbifold cohomology
algebra.

\section{$S^1$-Equivariant Floer Cohomology and Quantum Cohomology}
\label{sec:floer}

Floer cohomology should capture information about ``semi-infinite
cycles'' in the free loop space $L\stackX$.  Giving a rigorous
definition is not easy, particularly if one wants to define a theory
which applies beyond the toric setting, and we will not attempt to do
so here: various approaches to the problem can be found in
\cite{Iritani,Vla,Kapranov1, Kapranov2, CJS}.  Instead we will
indicate roughly how one might define Floer cohomology groups $\HF$ in
terms of Morse theory on a covering space of $\LX$, and explain how to
compute them.  We argue mainly by analogy with Morse theory on
finite-dimensional manifolds.  An excellent (and rigorous)
introduction to finite-dimensional Morse theory from a compatible
point of view can be found in \cite{AB}.  The material in this section
provides motivation and context for the rest of the paper, but most of
it is not rigorous mathematics: we do not discuss the topologies on
many of the spaces we consider, for example, and questions of
transversality and compactness are systematically ignored.  More
importantly, several key steps in the argument are plausible analogies
rather than rigorous proof.  None of the material in this section is
logically necessary, and so the reader may want to skip directly to
Section \ref{sec:j_fun}.

\subsection{Loops in $\stackX$}

Lupercio and Uribe have defined the loop groupoid of any topological
groupoid \cite{Lupercio--Uribe}.  As $\stackX$ can be represented by a
proper \'etale Lie groupoid \cite{Moerdijk}, this defines the
\emph{loop space} $L \stackX$.  Let $\cU = \CC^{n+1}- \{0\}$.  The
Lupercio--Uribe definition can be rephrased in the following
equivalent ways:
\begin{itemize}
\item[(A)] a loop in $\stackX$ is a pair $(\gamma,h)$ where
  $\gamma:[0,1] \to \cU$ is a continuous map and $h \in \Cstar$
  satisfies $\gamma(1) = h \gamma(0)$; loops $(\gamma_1,h_1)$ and
  $(\gamma_2,h_2)$ are isomorphic if there exists $k:[0,1] \to \Cstar$
  with $\gamma_2(x) = k(x) \gamma_1(x)$ for all $x \in [0,1]$ and $h_2
  = k(1) h_1 k(0)^{-1}$.
\item[(B)] a loop in $\stackX$ is a diagram 
  \begin{equation}
    \label{eq:descriptionB}
    \xymatrix{P \ar[r]^f \ar[d]& \cU \\ S^1}
  \end{equation}
  where $P \to S^1$ is a principal $\Cstar$-bundle and $f$ is a
  $\Cstar$-equivariant continuous map; an isomorphism between the
  loops 
  \begin{align*}
    \xymatrix{P_1 \ar[r]^{f_1} \ar[d]& \cU \\ S^1} && \text{and} &&
    \xymatrix{P_2 \ar[r]^{f_2} \ar[d]& \cU \\ S^1}    
  \end{align*}
  is an isomorphism $\phi:P_1 \to P_2$ of principal $\Cstar$-bundles
  such that the diagram
  \[
  \xymatrix{P_1 \ar[r]^{f_1} \ar[d] \ar[rd]^\phi & \cU \\ S^1 & P_2 \ar[u]_{f_2} \ar[l]}
  \]
  commutes.
\end{itemize}
The loop space $L\stackX$ can be thought of as an infinite-dimensional
K\"ahler orbifold, as follows.  A tangent vector to $L\stackX$ at
$(\gamma,h)$ is a vector field $v:[0,1] \to T\cU$ along $\gamma$ such
that $v(1) = h_\star v(0)$.  Weighted projective space is a K\"ahler
orbifold: let $\omega \in \Omega^2(\stackX)$ be the K\"ahler form on
$\stackX$ obtained by symplectic reduction from the standard K\"ahler
form on $\cU$ such that $\omega$ represents the class $P \in
H^2(\stackX;\CC)$, and let $g$ be the corresponding K\"ahler metric on
$\stackX$.  These structures induce a K\"ahler form on $L\stackX$:
\begin{align*}
  \Omega(u,v) = \int_0^1 \omega(u(t),v(t)) \, dt, && u,v \in
  T_{(\gamma,h)} L \stackX,
\end{align*}
and a Riemannian metric on $L\stackX$:
\begin{align*}
  G(u,v) = \int_0^1 g(u(t),v(t)) \, dt, && u,v \in
  T_{(\gamma,h)} L \stackX.
\end{align*}

\subsection{The Symplectic Action Functional} 
There is an $S^1$-action on $L\stackX$ given by rotation of loops (see
\cite{Lupercio--Uribe}).  This action is locally Hamiltonian with
respect to the K\"ahler form $\Omega$.  The moment map $m:L\stackX \to
S^1$ for this action, which is called the \emph{symplectic action
  functional}, is given as follows.  Every loop in $\stackX$ is the
boundary value of a representable continuous map $f:D \to \stackX$
from a possibly-stacky\footnote{Let $\Sigma$ be a Riemann surface,
  which may have a boundary.  By a ``possibly-stacky $\Sigma$'' we mean
  a reduced orbifold with coarse moduli space equal to $\Sigma$ and no
  stacky points on the boundary.} disc $D$.
The integral $\int_D f^\star \omega$ does not depend unambiguously on
the loop $\gamma$, because there are many possible choices of $D$ and
$f$, but the ambiguity in its value lies in the set
\[
\Pi = \bigg\{\int_S g^\star \omega \, \, \bigg| \, \, 
\begin{minipage}{2.5in}
  $S$ a possibly-stacky sphere,\\ $g:S \to \stackX$
  representable and continuous
\end{minipage}
\bigg\}.
\] 
Since $\RR/\Pi \cong S^1$, the map
\[
m: \gamma \longmapsto \int_D \gamma
\]
defines a circle-valued function on $L\stackX$.  This is the
symplectic action functional.  Pulling back the universal cover $\RR
\to S^1$ along the map $m:\LX \to S^1$ defines a covering $p:\tLX \to
\LX$ and a function $\mu: \tLX \to \RR$.  We can regard the covering
$\tLX$ as consisting of pairs $(\gamma, [D])$, where $\gamma$ is a
loop in $\stackX$ and $[D]$ is a relative homology class of
possibly-stacky discs $D$ with boundary $\gamma$.  The function $\mu$
gives the area of the disc $D$:
\[
\mu: (\gamma, [D]) \longmapsto \int_D \gamma.
\] 
We will study the Morse theory of $\mu$.

\begin{rem} \label{3rem:loopsbound} When applying this argument to
  other orbifolds $\cX$, one should consider only the subset of $L\cX$
  consisting of loops which bound possibly-stacky holomorphic discs.
  This condition does not arise here, as every loop in $\stackX$ is
  the boundary value of a representable continuous map $f:D \to
  \stackX$ from a disc $D$ with one possibly-stacky point at the
  origin.  To see this, observe that every loop in $\stackX$ is
  homotopic to a loop which lands entirely within the image of a
  co-ordinate chart
  \begin{align*}
    \big\{ [z_0:z_1:\cdots:z_n] \in \stackX \, \, \big| \, \, z_i = 1 \big\}
    && \text{for some $i$}
  \end{align*}
  and consequently (because these co-ordinate charts are contractible)
  that every loop in $\stackX$ is homotopic to a loop with image
  contained in one of the points
  \[
  \big\{ [z_0:z_1:\cdots:z_n] \in \stackX \, \, \big| \, \, \text{$z_j
    = 0$ for $j \ne i$, $z_i = 1$} \big\}.
  \]
  Such loops evidently bound representable continuous maps $f:D \to
  \stackX$, where $D$ is a disc with one possibly-stacky point at the
  origin, and the assertion follows.
\end{rem}
 
\subsection{Morse Theory}

As motivation, let us recall some key points from \cite{AB}.  Let
$(X,g)$ be a finite-dimensional Riemannian manifold and $f:X \to \RR$
a Morse--Bott function.  Let $X^{\text{cr}}_1,\ldots,X^{\text{cr}}_r$
denote the components of the critical set of $f$, $X^{\text{cr}} =
\coprod_i X^{\text{cr}}_i$, and let $\cM$ be the set of descending
gradient trajectories of $f$ (\emph{i.e.} of integral curves
$\gamma:\RR \to X$ for the vector field $-\mathop{grad}(f)$).  Under
reasonable conditions on $f$ and $g$, $\cM$ is a smooth
finite-dimensional manifold with a natural compactification $\clM$.  A
point of $\clM$ consists of a sequence of gradient trajectories
$\gamma_1(t),\ldots,\gamma_m(t)$, where $m \geq 1$, such that
\begin{align*}
\lim_{t \to \infty} \gamma_i(t) = \lim_{t \to -\infty}
\gamma_{i+1}(t), && 1 \leq i < m.
\end{align*}
There is an action of $\RR$ on $\cM$ by ``time translation'':
\begin{align*}
  \RR \times \cM & \longrightarrow \cM \\
  \big(s,\gamma(t)\big) & \longmapsto \gamma(s+t)
\end{align*}
and this extends to give an action on $\clM$.  Let $\gamma:\RR \to X$
be a descending gradient trajectory.  As $t \to \mp \infty$,
$\gamma(t)$ approaches critical points of $f$; this defines upper and
lower endpoint maps $u:\clM/\RR \to X^{\text{cr}}$ and $l:\clM/\RR \to
X^{\text{cr}}$. 

Chains in the Morse--Bott complex of $f$ are differential forms on the
critical set:
\[
C^{\text{MB}}_\bullet = \bigoplus_{i=1}^r
\Omega^\bullet(X^{\text{cr}}_i),
\]
where the grading on $\Omega^\bullet(X^{\text{cr}}_i)$ is shifted by
an integer which depends on the component $X^{\text{cr}}_i$ (see
\cite{AB}).  Consider the diagram
\[
\xymatrix{ X^{\text{cr}} & \clM/\RR \ar[r]^{l} \ar[l]_{u} &
  X^{\text{cr}}.}
\]
The differential in the Morse--Bott complex is the sum of the deRham
differential and a contribution from the space $\clM/\RR$ of gradient
trajectories:
\begin{align*}
d_{\text{MB}} \alpha = d_{\text{deRham}} \alpha + (-1)^j u_\star
l^\star \alpha && \text{for $\alpha \in \Omega^j(X^{\text{cr}})$.}
\end{align*}
The homology of the complex $(C_\bullet^{\text{MB}},d_{\text{MB}})$ is
the cohomology of $X$:
\begin{equation}
  \label{eq:MBiso}
  H_\bullet(C_\bullet^{\text{MB}},d_{\text{MB}}) \cong
  H^\bullet(X;\RR).
\end{equation}
Let $\alpha \in \Omega^\bullet(X^{\text{cr}}_i) \subset
C^{\text{MB}}_\bullet$, and let $A$ be a generic cycle in
$X^{\text{cr}}_i$ which is Poincar\'e-dual to $\alpha$.  Under the
isomorphism \eqref{eq:MBiso}, the class $[\alpha] \in
H_\bullet(C_\bullet^{\text{MB}},d_{\text{MB}})$ maps to the cohomology
class on $X$ which is Poincar\'e-dual to the cycle\footnote{More
  precisely, the Poincar\'e-dual cycle is the \emph{closure} of the
  locus 
  \[
  \bigcup \Big\{ \gamma(t) \,\, \Big | \,\, \text{$t \in \RR$, $\gamma$ is a
    gradient trajectory such that $l(\gamma) \in A$.}\Big\}.
  \]} swept out by gradient trajectories that end on $A$.  So, roughly
speaking, $\alpha \in \Omega^\bullet(X^{\text{cr}}_i)$ represents the
cohomology class dual to the cycle given by upward gradient flow from
$A \subset X^{\text{cr}}_i$.

Furthermore if $X$ is a finite-dimensional manifold with
$S^1$-action, $g$ is an $S^1$-invariant Riemannian metric on $X$, and
$f:X\to \RR$ is an $S^1$-invariant Morse--Bott function then, under
reasonable conditions on $f$ and $g$, we can compute the
$S^1$-equivariant cohomology of $X$ using the $S^1$-equivariant
Morse--Bott complex of $f$.  If we define chain groups
\[
C^{S^1,\text{MB}}_\bullet = \bigoplus_{i=1}^r
\Omega^\bullet(X^{\text{cr}}_i) \otimes H^\bullet_{S^1}(pt;\RR),
\]
and use the differential $d_{\text{MB}}$ as before, then
\[
H_\bullet(C_\bullet^{S^1,\text{MB}},d_{\text{MB}}) \cong
H^\bullet_{S^1}(X;\RR).
\]

\subsection{Floer Cohomology and $S^1$-Equivariant Floer Cohomology} 
Recall our setup
\[
\xymatrix{ {\tLX} \ar[r]^\mu \ar[d]_p & \RR \ar[d]^{t \mapsto \exp(2 \pi \sqrt{-1}
    t)} \\ {\LX} \ar[r]^m & S^1}
\]
where $m$ is the moment map for the $S^1$-action on $\LX$ given by
loop rotation.  We define the \emph{Floer cohomology} of $\LX$ to be
the homology of the Morse--Bott complex of $\mu$.  We will describe
the critical set of $\mu$ in a moment.  Gradient trajectories of
$\mu$, with respect to the induced K\"ahler metric $p^\star G$ on
$\tLX$, give paths of loops in $\stackX$ which sweep out holomorphic
cylinders.  It is this --- the link between Morse-theoretic gradient
trajectories and holomorphic curves --- which connects Floer
cohomology to Gromov--Witten theory.

The critical set of $\mu$ is a covering space of the critical set of
$m$.  As $m$ is a moment map, the critical set of $m$ coincides with
the $S^1$-fixed set on $\LX$.  This $S^1$-fixed set is canonically
isomorphic to the inertia stack $\inertia$ (see
\cite{Lupercio--Uribe}) and so the critical set of $\mu$ is a covering
space of $\inertia$.  The deck transformation group of this covering,
and of the covering $p:\tLX \to \LX$, is $\ZZ$: let $\CC[Q,Q^{-1}]$
denote the group ring of the group of deck transformations.  A deck
transformation changes the value of the function $\mu$ by an integer,
and we have
\[
\big( \text{critical set of $\mu$} \big) \cap \mu^{-1}(r) = 
\begin{cases}
  \text{a copy of $\PP(V^{\langle r \rangle})$} & \text{if $\langle r
    \rangle \in F$,} \\
  \varnothing & \text{otherwise}.
\end{cases}
\]
We will call a component of the critical set of $\mu$ which lies in
$\mu^{-1}(r)$ \emph{the component of the critical set at level $r$}.
A point in the component of the critical set at level $r$ is a pair
$(\gamma, [D])$ where $\gamma$ is an $S^1$-fixed loop in $\stackX$ and
$[D]$ is the homology class of a possibly-stacky disc bounding
$\gamma$ and having area $r$.  As $\gamma$ here is an $S^1$-fixed
loop, $[D]$ is in fact the homology class of a possibly-stacky
\emph{sphere}
 in $\stackX$ of area $r$.

The chain groups in the Morse--Bott complex for $\mu$ should be
\[
C_\bullet^{\text{MB}} = \Bigg(\bigoplus_{f \in F} Q^f
\Omega^\bullet\big(\PP\big(V^f\big)\big) \Bigg) \otimes \CC[Q,Q^{-1}].
\]
Here we introduced fractional\footnote{These fractional shifts will
  play an essential role later --- see \eqref{approximation} and the
  discussion thereafter.}  powers $Q^f$ so that an element $\alpha Q^r
\in C_\bullet^{\text{MB}}$, where $\alpha \in
\Omega^\bullet\big(\PP\big(V^f\big)\big)$ and $r \in \QQ$, is a
differential form $\alpha$ on the component of the critical set at
level $r$.  The grading on the chain groups is defined by
\[
\deg \big(\alpha Q^r\big) = \deg \alpha + \age \PP\big(V^f\big) + (w_0 + \ldots +
w_n) r.
\]
Note that $\deg (\alpha Q^r) \in \ZZ$.  As before, the
differential in the Morse--Bott complex should be given by
\[
d_{\text{MB}} \theta = d_{\text{deRham}} \theta \pm u_\star
l^\star \theta,
\]
where $u$ and $l$ fit into the diagram
\[
\xymatrix{  X^{\text{cr}} & \clM/\RR \ar[r]^{l} \ar[l]_{u} & X^{\text{cr}}  .}
\]
In this case the space $\cM$ of descending gradient trajectories, each
of which gives a holomorphic map $\cC \to \stackX$ from a cylinder
$\cC$, admits an $S^1$-action coming from the reparametrization of
$\cC$.  This $S^1$-action extends to an $S^1$-action on $\clM$, which
commutes with the $\RR$-action on $\clM$.
\begin{equation}
  \label{eq:correctionsdie}
  \xymatrix{X^{\text{cr}} & \clM/\RR \ar[r]^{l} \ar[l]_{u} \ar[d] & X^{\text{cr}}\\ 
    & \clM/(\RR \times S^1) \ar[lu] & }
\end{equation}
The upper and lower endpoint maps $u$, $l$ are $S^1$-equivariant, and
so for each $\theta \in \Omega^\bullet(X^{\text{cr}})$ we have
$u_\star l^\star \theta = 0$: we can compute the pushforward along $u$
by first pushing forward along the vertical map in
\eqref{eq:correctionsdie}, and this pushforward sends the
$S^1$-invariant differential form $l^\star \theta$ to zero.  Thus in
this case we should have $d_{\text{MB}} \theta = d_{\text{deRham}}
\theta$, and so
\[
H_\bullet(C_\bullet^{\text{MB}},d_{\text{MB}})\\
= \Bigg(\bigoplus_{f \in F} Q^f
H^\bullet\big(\PP\big(V^f\big);\RR\big) \Bigg) \otimes \CC[Q,Q^{-1}]
\]
as graded vector spaces.  Here the grading on
$H^\bullet\big(\PP\big(V^f\big)\big)$ is shifted by the age of
$\PP(V^f)$, and the degree of $Q$ is $w_0 + \cdots + w_n$.  

It follows, as indicated in Remark~\ref{1rem:HF}, that after
completing the group ring $\CC[Q,Q^{-1}]$ we can identify $\HF =
H_\bullet(C_\bullet^{\text{MB}},d_{\text{MB}})$ with the free
$\NovZ$-submodule of $\HorbLinv$ with basis \eqref{1eq:basis}.  The
$\NovZ$-module structure here arises from the action of deck
transformations on $\tLX$.  Let $z$ be the first Chern class of the
tautological line bundle over $BS^1$, so that $H^\bullet_{S^1}(pt) =
\CC[z]$.  Identical arguments and conventions suggest that the
$S^1$-equivariant Floer cohomology $\eqHF$ should be the free
$\NovzZ$-submodule of
\[
\Horb \otimes
\CC[z][\![Q^{1/\mathrm{lcm}(w_0,\ldots,w_n)}]\!][Q^{-1}]
\]
with basis \eqref{1eq:basis}.  

\subsection{Floer Cohomology and Small Quantum Cohomology}

We think of elements of $\HF$ as representing semi-infinite cycles in
$\tLX$, as follows. Recall that gradient trajectories of $\mu:\tLX \to
\RR$ sweep out holomorphic cylinders in $\stackX$.  Recall further
that we are using bases $\phi_1,\ldots,\phi_N$ and
$\phi^1,\ldots,\phi^N$ for $\Horb$ such that
$\(\phi^i,\phi_j\)_{\text{orb}} = \delta^i_{\phantom{i}j}$.  Suppose
that $\phi_\beta \in \Horb$ is supported on $\PP(V^f) \subset
\inertia$, and let $B$ be a generic cycle in $\PP(V^f)$ which is
Poincar\'e-dual to $\phi_\beta$.  The Floer cohomology class
$\phi_\beta Q^r \in \HF$ represents the semi-infinite cycle in $\tLX$
swept out by upward gradient flow from the copy of $B$ in the
component\footnote{Note that $\langle r \rangle = f$, and so the
  component of the critical set at level $r$ is a copy of $\PP(V^f)$.}
of the critical set at level $r$.  The projection of this
semi-infinite cycle to $\LX$ consists of loops\footnote{More
  precisely, the projection consists of the closure of the set of such
  loops.  In the rest of this section, we will ignore such
  distinctions.} in $\stackX$ which bound a holomorphic disc $\{|z|
\leq 1\} \to \stackX$ with a possibly-stacky point at the origin such
that the $S^1$-fixed loop defined by the origin of the disc lies in $B
\subset \inertia$.

From this point of view, it is not obvious that $\HF$ should carry a
ring structure: the transverse intersection of two semi-infinite
cycles need not be semi-infinite, so we should not expect an
intersection product here.  But the transverse intersection of a
finite-codimension cycle with a semi-infinite cycle will be
semi-infinite, and this should give a map
\[
H^\bullet(\tLX) \otimes \HF \longrightarrow \HF.
\]
Evaluation at $1 \in S^1$ gives a map $\tLX \to \stackX$, and via
pull-back we get a map
\begin{align}
H^\bullet(\stackX;\CC) &\otimes \HF  \longrightarrow \HF \label{modulemap}\\
\phi_\alpha &\otimes \phi_\beta Q^r  \longmapsto \sum_{d \in \QQ}
\sum_\gamma n(d)_{\alpha \beta}^{\phantom{\alpha \beta} \gamma}
\phi_\gamma Q^{d + r} \notag
\end{align}
which\footnote{Note that \eqref{modulemap} involves the subspace
  $H^\bullet(\stackX;\CC) \subset \Horb$ and not the full orbifold
  cohomology group $\Horb$.} commutes with the action of $\NovZ$.
The structure constants of this map have a geometric interpretation,
as follows.  If everything intersects transversely, the structure
constant $n(d)_{\alpha \beta}^{\phantom{\alpha \beta} \gamma}$ should
count the number of isolated points in the intersection of three
cycles in $\tLX$:
\begin{itemize}
\item[(a)] the finite-codimension cycle corresponding to $\phi_\alpha$;
\item[(b)] the semi-infinite cycle corresponding to $\phi_\beta Q^r$;
\item[(c)] a semi-infinite cycle representing the element of Floer
  {\em homology} corresponding to $\phi_\gamma Q^{d + r}$.
\end{itemize}
Cycle (a) is the pre-image in $\tLX$ of the cycle in $\LX$ consisting
of loops such that the point $1 \in S^1$ maps to a generic cycle in
$\stackX$ Poincar\'e-dual to $\phi_\alpha$.  Cycle (b) was described
above.  Cycle (c) is swept out by {\em downward} gradient flow from an
appropriate cycle in the component of the critical set at level $d+r$.
Its projection to $\LX$ consists of loops which bound a holomorphic
disc $\{|z| \geq 1\} \to \stackX$ with a possibly-stacky point at
$\infty$ such that the $S^1$-fixed loop defined by the point $\infty$
lies in a generic cycle in $\inertia$ Poincar\'e-dual to
$\phi^\gamma$.  So $n(d)_{\alpha \beta}^{\phantom{\alpha \beta}
  \gamma}$ counts --- or, in the non-transverse situation, gives a
virtual count of --- the number of isolated holomorphic spheres in
$\stackX$ of degree $d \in \QQ$ carrying exactly two possibly-stacky
points $\{0,\infty\}$ and incident at the points $\{0,1,\infty\}$ to
generic cycles in $\inertia$ Poincar\'e-dual respectively to
$\phi_\beta$, $\phi_{\alpha}$, and $\phi^\gamma$.  In other words, the
structure constants $n(d)_{\alpha \beta}^{\phantom{\alpha \beta}
  \gamma}$ of the map \eqref{modulemap} coincide with the structure
constants \eqref{2eq:smallQC} of the small orbifold quantum cohomology
algebra.

\begin{rem}
  This shows that small quantum orbifold multiplication by a class in
  the untwisted sector $H^\bullet(\stackX;\CC) \subset \Horb$ can be
  thought of as an operation on Floer cohomology.  It would be
  interesting to find an interpretation of multiplication by other
  orbifold cohomology classes in these terms.
\end{rem}

\subsection{The $\cD$-Module Structure on $S^1$-Equivariant Floer
  Cohomology}

In this section we explain why the $S^1$-equivariant Floer cohomology
$\eqHF$ should carry a natural $\cD$-module structure.  Recall that
$\Omega$ is the K\"ahler form on $\LX$ induced by the K\"ahler
structure on $\stackX$, and that we consider the covering space
$p:\tLX \to \LX$.  We have $[\Omega] = \ev_1^\star P$.  The form
$p^\star \Omega$ is not equivariantly closed, so it does not define an
$S^1$-equivariant cohomology class on $\tLX$, but $p^\star \Omega + z
\mu$ is equivariantly closed --- this follows from the fact that $m$
is a moment map.  Let $\wp$ be the class of $p^\star \Omega + z \mu$
in $H^2_{S^1}(\tLX)$.  Consider the map $\bP:\eqHF \to \eqHF$ given by
multiplication by $\wp$, and the map $\bQ:\eqHF \to \eqHF$ given by
pull-back by the deck transformation $Q^{-1}$.  Since
\[
(Q^{-1})^\star \wp = \wp - z
\]
we have $[\bP,\bQ] = z\bQ$.  In other words, if we define $\cD$ to be
the Heisenberg algebra
\begin{align*}
  \cD = \CC[z][\![\bQ]\!][\bQ^{-1}] \langle \bP \rangle && \text{such
    that} && [\bP,\bQ] = z\bQ,
\end{align*}
then $\eqHF$ should carry the structure of a $\cD$-module where $\bQ$
acts by pull-back by $Q^{-1}$ and $\bP$ acts by multiplication by
$\wp$.

In the non-equivariant limit ($z \to 0$) this structure degenerates to
a $\NovZ[P]$-module structure on $\HF$, where $P$ acts via
\eqref{modulemap}.  Thus we can recover the part of the small orbifold
quantum cohomology algebra generated by $P$ --- which, as we will see
below, is the whole thing --- from the $\cD$-module structure on
$\eqHF$.  It is clear that $\HF$ should be generated as a
$\NovZ[P]$-module by $\{Q^f \fun_f\}$, so we expect $\eqHF$ to be
finitely generated as a $\cD$-module.  Our analysis below will show
that $\eqHF$ is of rank one, generated by $\fun_0 Q^0$.  This
generator is Givental's ``fundamental Floer cycle'' --- it represents
the semi-infinite cycle in $\tLX$ swept out by upward gradient flow
from the component of the critical set at level $0$.  The projection
to $\LX$ of the fundamental Floer cycle consists of all loops which
bound holomorphic discs with a possibly-stacky point at the origin.

The link between Floer cohomology and Gromov--Witten theory appears
here as a conjectural $\cD$-module isomorphism between $\eqHF$ and the
$\cD$-module generated by the small $J$-function.  We have seen how
$\cD$ acts on $\eqHF$.  Another realization of $\cD$ is by
differential operators
\begin{align*}
  \bP \colon f & \mapsto z {\partial f \over \partial t}  & \bQ
  \colon f & \mapsto Q e^t f
\end{align*}
acting on the space of analytic functions $f:\CC \to \HorbLinv \otimes
\CC(\!(z^{-1})\!)$.  The small $J$-function is such a function (see
Section~\ref{2sec:smallJ}) and so it generates a $\cD$-module;
relations in this $\cD$-module are differential equations satisfied by
$J_{\stackX}(t)$ (see equations \eqref{2eq:diffeq},
\eqref{2eq:smallQDEs} and the discussions thereafter).  We will make
use of this conjectural $\cD$-module isomorphism in the next Section,
where we write down a concrete model for $\eqHF$ as a $\cD$-module and
then identify the fundamental Floer cycle in this model with the small
$J$-function $J_{\stackX}(t)$.  This will give a conjectural formula
for the small $J$-function.

\subsection{Computing the $\cD$-Module Structure}

As we lack a concrete model for $\tLX$, we consider instead the space
of \emph{polynomial loops}
\[
\Lpol = \left\{ \begin{pmatrix} f^0,\ldots,f^n \end{pmatrix} \mid
  \mbox{$f^i \in \CC[t,t^{-1}]$, not all the $f^i$ are zero} \right\}
\slash \Cstar
\] 
where $\alpha \in \Cstar$ acts on a vector-valued Laurent polynomial
as:
\[
\begin{pmatrix} f^0,\ldots,f^n \end{pmatrix} \longmapsto
\begin{pmatrix} \alpha^{-w_0} f^0,\ldots,\alpha^{-w_n} f^n \end{pmatrix}.
\]
The space $\Lpol$ is quite different from $\tLX$ --- it is, for
example, certainly not a covering space\footnote{The ``obvious map''
  $\Lpol \to \LX$, given by restricting a polynomial map $f(t)$ to the
  circle $\{t \in \CC : |t|=1\}$ and filling in where necessary using
  continuity, is not even continuous.} of $\LX$.  But $\Lpol$ is in
some ways a good analog for $\tLX$.  We will see below that there is
an $S^1$-action on $\Lpol$ such that the $S^1$-fixed subset is a
covering space of the inertia stack $\inertia$ with deck
transformation group $\ZZ$.  So for computations involving quantities
which localize to the $S^1$-fixed set --- such as $S^1$-equivariant
semi-infinite cohomology --- $\Lpol$ is a good substitute for $\tLX$.
Working by analogy with the discussion in the previous Section, we now
construct an action of $\cD$ on the ``$S^1$-equivariant semi-infinite
cohomology'' of $\Lpol$.  This will be our concrete model for $\eqHF$.

The space $\Lpol$ is an infinite-dimensional weighted projective
space.  It carries an $S^1$-action coming from loop rotation, which is
Hamiltonian with respect to the Fubini-Study form $\Omega' \in
\Omega^2(\Lpol)$.  The moment map for this action is
\[
\mu' \colon \left[ 
\begin{pmatrix} 
\sum_{k\in\ZZ} a^0_k t^k,&\ldots, &\sum_{k\in\ZZ} a^n_k t^k
\end{pmatrix}
\right]
\longmapsto
- {\sum_{0 \leq l \leq n} \sum_{k\in\ZZ} k |a^l_k|^2 \over \sum_{0 \leq
    l \leq n} \sum_{k\in\ZZ} w_l |a^l_k|^2 }.
\]
A polynomial loop
\[
\left[
\begin{pmatrix}
  f^0(t),& \ldots, &f^n(t)
\end{pmatrix}
\right]
\in \Lpol
\]
is fixed by loop rotation if and only if
\[
\begin{pmatrix}
  f^0(\lambda t),& \ldots, &f^n(\lambda t)
\end{pmatrix}
=
\begin{pmatrix}
  \alpha(\lambda)^{-w_0} f^0(t),& \ldots, &\alpha(\lambda)^{-w_n} f^n(t)
\end{pmatrix}
\]
for all $\lambda\in S^1$ and some possibly multi-valued function
$\alpha(\lambda)$.  We need $\alpha(\lambda)=\lambda^{-k/w_i}$ for some
integer $k$, so components of the $S^1$-fixed set are indexed by
\[
\widetilde{F} = \left\{ {k \over w_i} \, \bigg\vert \, \, k \in \ZZ, 0 \leq i \leq n \right\}.
\]
For $r \in \widetilde{F}$, the corresponding $S^1$-fixed component
\[
\Fix_r = \Big\{
\left[
\begin{pmatrix}
  b_0 t^{w_0 r}, & \ldots, & b_n t^{w_n r}
\end{pmatrix}
\right] \in \Lpol \,\, \Big| \,\,
\mbox{$b_i = 0$ unless $w_i r \in \ZZ$}
\Big\}
\]
is a copy of the component $\PP(V^{\langle r \rangle})$ of the inertia
stack, and the value of $\mu'$ on this fixed component is $-r$.  The
normal bundle to $\Fix_r$ is
\[
\bigoplus_{i=0}^{i=n} \bigoplus_{\substack{j \in \ZZ \\ j \ne w_i r}}
\cO\big(w_i P + (j - w_i r) z\big),
\]
where $\cO(a P + b z)$ denotes the bundle $\cO(a)$ on $\Fix_r =
\PP(V^{\langle r \rangle})$ which has weight $b$ with respect to loop
rotation.

Let $\wp'$ be the class of $\Omega' + z \mu'$ in $H^2_{S^1}(\Lpol)$,
so that
\[
H^\bullet_{S^1}(\Lpol) = \CC[z,\wp'],
\]
and introduce an action of $\ZZ$ on $\Lpol$ by ``deck
transformations'':
\begin{multline*}
Q^m  \colon  \left[ 
\begin{pmatrix} 
\sum_{k\in\ZZ} a^0_k t^k,&\ldots, &\sum_{k\in\ZZ} a^n_k t^k
\end{pmatrix}
\right]
\longmapsto \\ 
\left[ 
\begin{pmatrix} 
\sum_{k\in\ZZ} a^0_k t^{k-m w_0},&\ldots, &\sum_{k\in\ZZ} a^n_k t^{k -
  m w_n}
\end{pmatrix}
\right],
\quad 
m \in \ZZ.
\end{multline*}
The deck transformation $Q^m$ changes the value of $\mu'$ by $m$, and
sends $\Fix_r$ to $\Fix_{r-m}$.  We let $\bQ$ act on
$H^\bullet_{S^1}(\Lpol)$ by pull-back by $Q^{-1}$, and $\bP$ act on
$H^\bullet_{S^1}(\Lpol)$ by cup product with $\wp'$.  As
\[
(Q^{-1})^\star \wp' = \wp' - z,
\]
so that $[\bP,\bQ] = z \bQ$, this gives an action of $\cD$ on
$H^\bullet_{S^1}(\Lpol)$.

We now consider the ``$S^1$-equivariant semi-infinite cohomology'' of
$\Lpol$.  We will work formally, representing semi-infinite cohomology
classes by infinite products in $H^\bullet_{S^1}(\Lpol)$.  These
products, interpreted na\"{\i}vely, definitely diverge, but one can
make rigorous sense of them by considering them as the limits of
finite products and at the same time considering $\Lpol$ as the limit
of spaces of Laurent polynomials of bounded degree.  This is explained
in \cite{Givental:homological,Iritani}.  Recall that the fundamental
Floer cycle in $\tLX$ consists (roughly speaking) of loops which bound
holomorphic discs.  The analog of the fundamental Floer cycle in
$\Lpol$ is the cycle of Laurent polynomials which are regular at
$t=\infty$.  We represent this by the infinite product
\[
\Delta = \prod_{i=0}^{i=n} \prod_{k>0} (w_i \wp' + k z).
\]
To interpret this, observe that the Fourier coefficient $a^i_k$ of the loop
\[
\left[ 
  \begin{pmatrix}
    \sum_{k \in \ZZ} a^0_k t^k, & \ldots, &     \sum_{k \in \ZZ} a^n_k t^k 
  \end{pmatrix}
\right]
\in \Lpol
\]
gives a section of the bundle $\cO(w_i)$ over 
\[
\Lpol \iso
\PP(\ldots,w_n,w_0,w_1,\ldots,w_n,w_0,w_1,\ldots,w_n,w_0,\ldots)
\]
which has weight $k$ with respect to loop rotation.  Our candidate for
the Floer fundamental cycle is cut out by the vanishing of the
$a^i_k$, $k>0$, and so $\Delta$ is a candidate for the
$S^1$-equivariant Thom class of its normal bundle --- that is, for its
$S^1$-equivariant Poincar\'e-dual.  We have
\begin{equation}
  \label{3eq:Dmodulerelations}
  \prod_{i=0}^{i=n} \prod_{j=0}^{j=w_i-1} \( w_i \bP - j z \) \, \Delta
  = \bQ \, \Delta.
\end{equation}
This is an equation in the $S^1$-equivariant semi-infinite cohomology
of $\Lpol$, regarded as a $\cD$-module via the actions of $\bP$ and
$\bQ$ defined above.  As a $\cD$-module, the $S^1$-equivariant
semi-infinite cohomology of $\Lpol$ is generated by $\Delta$.

We cannot directly identify $\Delta$ with the small $J$-function, as
the $\cD$-module generated by $\Delta$ involves shift operators
\begin{align*}
\bP \colon  g(\wp') & \mapsto \wp' g(\wp') & \bQ \colon g(\wp') & \mapsto g(\wp'-z)
\end{align*}
whereas that generated by the small $J$-function involves differential
operators
\begin{align*}
  \bP \colon  f(t) & \mapsto z {\partial f \over \partial t}  &
  \bQ \colon  f(t) & \mapsto Q e^t f(t).
\end{align*}
We move between the two via a sort of Fourier transform.  We expect,
by analogy with the Atiyah--Bott localization theorem
\cite{Atiyah--Bott}, that there should be a localization map $\Loc$
from localized $S^1$-equivariant semi-infinite cohomology of $\Lpol$
to the cohomology $H^\bullet_{S_1}(\Lpol^{S^1})\otimes \CC(z)$ of the
$S^1$-fixed set.  We consider
\begin{equation} \label{loc}
\Loc\(e^{\wp' t/z} \Delta\)
\end{equation}
as this should satisfy
\begin{align*}
  \bP \Loc\(e^{\wp' t/z} \Delta\) &= z {\partial \over \partial t}
  \Loc\(e^{\wp' t/z} \Delta\) \\
  &= \Loc\(e^{\wp' t/z} \, \wp' \Delta\) \\
  &= \Loc\(e^{\wp' t/z} \, \bP  \Delta\) \\
  \bQ  \Loc\(e^{\wp' t/z} \Delta\) &= Q e^t \Loc\(  e^{\wp' t/z} \Delta\) \\
  &= e^t \Loc\((Q^{-1})^\star  \(e^{\wp' t/z} \Delta\)\) \\
  &= \Loc\(e^{\wp' t/z} \, \bQ \Delta\).
\end{align*}
The class $\wp'\in H^2_{S^1}(\Lpol)$ restricts to the class
$c_1(\cO(1)) - zr \in H^2(\Fix_r)$, and we can write this as the
Chen--Ruan orbifold cup product
\[
(P-zr)\,\fun_{\langle -r \rangle}.
\]
Thus $\Loc(e^{\wp' t/z}\Delta)$ should be something like
\begin{equation}\label{approximation}
\sum_{r\in\widetilde{F}} Q^{-r} e^{Pt/z} e^{-rt}
{\prod_{i=0}^{i=n} \prod_{k>0} \( w_i P + (k-w_i r) z \) \over
\prod_{i=0}^{i=n} \prod_{\substack{j \in \ZZ \\ j \ne w_i r}} \( w_i P
+ (k-w_i r)z \)} \fun_{\langle -r \rangle} 
\end{equation}
where the numerator records the restriction of $\Delta$ to $\Fix_r$
and the denominator stands for the $S^1$-equivariant Euler class of
the normal bundle to $\Fix_r$.  We need to make sense of this
expression.

Note first that if $r>0$, the numerator in \eqref{approximation} is
divisible by $P^{\dim_{\langle r \rangle} +1}$ and hence vanishes for dimensional
reasons.  So our expression is
\[
\sum_{\substack{r\in\widetilde{F}\\r \geq 0}} Q^r e^{P t/z} e^{r t}
\prod_{i=0}^n 
{1 \over \prod_{\substack{b \colon \langle b \rangle = \langle r w_i
      \rangle \\ 0 < b \leq w_i r}} \( w_i P + b z\)}
{1 \over \prod_{\substack{b \colon \langle b \rangle = \langle r w_i
      \rangle \\ b<0}} \( w_i P + b z\)}
\fun_{\langle r \rangle}.
\]
This expression still does not make sense due to the divergent
infinite product on the right.  We ``regularize'' it by simply
dropping these factors --- which depend on $r$ only through $\langle r
\rangle$ --- and multiplying by $z$, obtaining the {\em $I$-function}:
\[
I(t)= z \, e^{Pt/z} 
\sum_{\substack{r\in\widetilde{F}\\r \geq 0}} Q^r e^{r t}
{1 \over \prod_{i=0}^n \prod_{\substack{b \colon \langle b \rangle = \langle r w_i
      \rangle \\ 0 < b \leq w_i r}} \( w_i P + b z\)}
\fun_{\langle r \rangle}.
\]
This is a formal function of $t$ taking values in $\HorbL$.  It satisfies
\[
\prod_{i=0}^{i=n} \prod_{j=0}^{j=w_i-1} \( w_i z {\partial \over
  \partial t} - j z \) \, I
= Q e^t \, I,
\]
so the $\cD$-modules generated by $\Delta$ and by $I$ are isomorphic
(see \eqref{3eq:Dmodulerelations}).  We conjecture that this
$\cD$-module is isomorphic to the $\cD$-module generated by the small
$J$-function, and that
\[
J_{\stackX}(t) = I(t).
\]

\section{Calculation of the Small $J$-Function}
\label{sec:j_fun}

\subsection{Summary: the Basic Diagram}
\label{sec:summ-basic-diagr}

In this section we describe a certain commutative diagram of stacks
with $\Cstar$-action which lies at the heart of our proof of
Theorem~\ref{1thm:smallJ}.  We begin by showing that for each
genus-zero one-pointed twisted stable map to $\stackX$, the component
of $\inertia$ to which the marked point maps is determined by the
degree of the map.

\begin{lem}
  \label{4lem:7} Fix a positive rational number $d>0$.
  \begin{enumerate}
  \item If the moduli stack $\p^\mathbf{w}_{0,0,d}$ is nonempty, then
    $d$ is an integer.
  \item If the moduli stack $\(\p^{\mathbf{w}}\)^f_{0,1,d}$ is nonempty,
    then $f=\langle -d\rangle$.
  \end{enumerate}
\end{lem}

\begin{proof}
  Let $\cC$ be a balanced twisted curve, and assume that there is a
  stable representable morphism $\varphi \colon \cC \to \p^\mathbf{w}$
  of degree $d$:
  \[
  \int_{\cC} \varphi^\star \o(1)=d.
  \]
  Applying Riemann--Roch for twisted curves
  \cite{AGV:2}*{Theorem~7.2.1}, we find that
  \[
  \chi(\cC,\varphi^\star \cO(1)) = 
  \begin{cases}
    1 + d & \text{in case (1)} \\
    1 + d - \langle -f \rangle & \text{in case (2)}.
  \end{cases}
  \]
  As $\chi(\cC,\varphi^\star \cO(1))$ is an integer, the result follows.
\end{proof}

\begin{nt}
  The lemma says that in $\(\stackX\)^f_{0,1,d}$ we always have
  $f=\langle -d\rangle$.  It is therefore safe to drop $f$ from the
  notation, and we do so in what follows.  Fix now $d=m/r$ in lowest
  terms and write $f = \langle -d\rangle \in F$.  We introduce the
  following notation:
  
  \begin{enumerate}
  \item $M_d=\Mbar_{0,1} (\p^\mathbf{w}, d)$ is, using the notation of
    \cite{AGV:2}, the moduli stack of genus-zero one-pointed balanced
    twisted stable morphisms of degree $d$ to $\p^\mathbf{w}$
    \emph{with section to the gerbe marking}.  There are maps
    \[
    \xymatrix{\mathcal{U} \ar[rr]^\varphi \ar[d]^\pi & & \p^\mathbf{w}\\
      M_d\ar@/^1pc/[u]^{\sigma}\ar[rr]_{\ev_1} & & \p (V^f)\ar[u]}
    \] 
    where $\pi\colon \mathcal{U}\to M_d$ denotes the universal family,
    $\sigma \colon M_d \to \mathcal{U}$ the section, and $\ev_1 \colon
    M_d \to \p(V^f)$ the evaluation map. As usual, we write $\psi_1 =
    c_1(L_1)$ where $L_1$ is the universal cotangent line at the
    marked point.
  \item $G_d$ is the \emph{graph space} of degree $d$; its definition
    depends on whether or not $d$ has a nontrivial fractional part:
    \[
    G_d=
    \begin{cases}
      \Mbar_{0,1} \bigl(\p^\mathbf{w}\times \p^{1,r}, \,d\times
      \frac1{r} \bigr) \quad & \text{if}\; \langle d \rangle >0\\
      \Mbar_{0,0} \bigl(\p^\mathbf{w}\times \p^1,\, d\times 1\bigr)
      \quad & \text{if}\; \langle d \rangle = 0
    \end{cases}
    \]
    More precisely, if $\langle d \rangle > 0$ then $G_d$ denotes the
    moduli stack of graphs with the following specified
    \emph{character} at the marked point: a point of $G_d$ is a pair
    of morphisms $(f_1,f_2)\colon C \to \p^\mathbf{w}\times \p^{1,r}$
    where $f_1\colon C \to \p^\mathbf{w}$, $f_2 \colon C \to
    \p^{1,r}$, and we require that $f_1$ evaluates in $\p (V^{\langle
      - d\rangle}) \subset \cI \stackX$ and $f_2$ evaluates in $\p
    (V^{\frac{r-1}{r}}) \subset \cI \PP^{1,r}$.  In other words,
    denoting by $x \in C$ the marked point,
    \begin{equation}
      \label{4eq:9}
      \begin{split}
        \Aut_C(x) & \longrightarrow 
        \Aut_{\stackX}\(f_1(c)\)\times \Aut_{\PP^{1,r}}(0) \\
        e^{\frac{2\pi \sqrt{-1}}{r}} & \longmapsto 
        \(e^{2 \pi \sqrt{-1} f}, e^{-\frac{2\pi \sqrt{-1}}{r}}\).
      \end{split}
    \end{equation}
    As a result of this choice, the marked point $x\in C$ is
    constrained to lie above the orbifold point $0\in \p^{1,r}$.  Note
    again that, if $\langle d \rangle>0$, our graphs have a gerbe
    marking and $G_d$ is a moduli stack of morphisms with section to
    the gerbe marking.
  \item $L_d$ is the stack of polynomial morphisms $\p^{1, r}\to
    \p^\mathbf{w}$ of degree $d$.  This is described in detail in
    Section~\ref{4sec:space-polyn-maps}.
  \end{enumerate}
\end{nt}

\begin{nt} In what follows 
  \begin{enumerate}
  \item all group actions are strict (see \emph{e.g.} \cite{Romagny});
  \item all stacks which we consider are Deligne--Mumford stacks,
    except where we explicitly say ``Artin stack'';
  \item we write ``stable morphism'' instead of ``balanced twisted
    stable morphism'';
  \item by ``part'' we mean ``union of connected components''.
  \end{enumerate}
\end{nt}

The action of the group $\Cstar$ on $\CC^2$
\begin{align} \label{4eq:action}
  \lambda \colon (s_0,s_1) &\longmapsto (s_0, \lambda^{-1} s_1) && 
  \lambda \in \Cstar
\end{align}
descends to give an action of $\Cstar$ on $\PP^{1,r} = \PP(1,r)$.
This action induces actions on the stacks $G_d$ and $L_d$; see below
for additional details and discussion.  

\begin{nt}
  Let $z$ be the first Chern class of the tautological line bundle on
  $B \Cstar$, so that $H_{\Cstar}^\bullet(\{\text{pt}\}) = \CC[z]$.
\end{nt}

\begin{thm}
  \label{4thm:1}
  There is a commutative diagram of stacks with
  $\C^\times$-action: 
\[
\xymatrix{
G_d \ar[r]^u & L_d \\
M_d \ar[r]_{\ev_1}\ar[u]^\iota & \p(V^f)\ar[u]_j\\ 
}
\]
such that the following properties hold:
\begin{enumerate}
\item The inclusion $j\colon \p(V^f) \hookrightarrow L_d$ is a
  connected component of the $\C^\times$-fixed substack, and the
  $\C^\times$-equivariant Euler class of its normal bundle is
  \[
  e(N_j)= \prod_{i=0}^n \prod_{\substack{b \colon \langle b \rangle =
      \langle dw_i \rangle \\0<b\leq dw_i}}(w_iP+bz).
  \]
\item The inclusion $\iota \colon M_d \hookrightarrow G_d$ is the part
  of the $\C^\times$-fixed substack of $G_d$ mapping to $\p(V^f)$.
  The canonical perfect obstruction theory on $M_d$ coincides with the
  perfect obstruction theory inherited from $G_d$, and the
  $\C^\times$-equivariant Euler class of the virtual normal bundle to
  $M_d$ is:
  \[
  e(N_\iota^\text{vir})=z(z-\psi_1).
  \]
\item The morphism $u$ is ``virtually birational''; in other words,
  when $G_d$ is endowed with its canonical perfect obstruction theory
  and $L_d$ with its intrinsic perfect obstruction theory, then
\[ 
u_\star \mathbf{1}_{G_d}^\text{vir}= \mathbf{1}_{L_d}.
\]
\end{enumerate}
\end{thm}

More details on obstruction theory can be found below.

\begin{cor}
  \label{4cor:1}
  Theorem~\ref{1thm:smallJ} of the Introduction holds.  That is, we have
  the following formula for the small $J$-function of $\p^\mathbf{w}$:
  \[
  J_{\stackX}(t)
  = z \, e^{P t/z}\sum_{\substack{d \colon d \geq 0 \\ \langle d \rangle \in F}}
  \frac{Q^d e^{d t} \mathbf{1}_{\langle d \rangle}}{
    \prod_{i=0}^n\prod_{\substack{b \colon \langle b \rangle = \langle
        dw_i \rangle \\ 0<b\leq dw_i}} (w_iP+bz)}.
  \]
\end{cor}

\begin{proof}[Proof of the Corollary]
We calculate using the basic diagram and properties stated in
Theorem~\ref{4thm:1}:
\begin{align}
  \mathbf{1}_f &= j^\star \mathbf{1}_{L_d}\notag \\
  & =j^\star u_\star
  \mathbf{1}_{G_d}^\text{vir}  \label{4eq:line_1}\\
  &=j^\star u_\star \iota_\star \biggl(\mathbf{1}_{M_d}^\text{vir}\cap
  \frac{1}{e(N_\iota^\text{vir})} \biggr) \label{4eq:line_2} \\ \notag
  &= j^\star u_\star \iota_\star \biggl(\mathbf{1}_{M_d}^\text{vir}\cap
  \frac{1}{z(z-\psi_1)} \biggr) \\ \notag &= j^\star j_\star
  \(\ev_1\)_\star \biggl(\mathbf{1}_{M_d}^\text{vir}\cap
  \frac{1}{z(z-\psi_1)} \biggr) \\ \notag &=e(N_j) \(\ev_1\)_\star
  \biggl( \mathbf{1}_{M_d}^\text{vir} \cap
  \frac1{z(z-\psi_1)}\biggr). \label{4eq_line3}
\end{align}
Equation~\eqref{4eq:line_1} here holds because $u$ is virtually
birational. Equation~\eqref{4eq:line_2} follows from the virtual
localization formula of Graber and Pandharipande
\cite{Graber--Pandharipande} and the fact that $M_d$ is the part of
the $\C^\times$-fixed substack of $G_d$ which maps to $\p(V^f)$. The
Graber--Pandharipande formula requires all stacks to admit a global
equivariant embedding in a smooth stack; the main result of
\cite{AGOT} shows that this is true here. From this, we conclude that:
\[
(I\circ \ev_1)_\star \biggl(\mathbf{1}_{M_d}^\text{vir} \cap
\frac1{z(z-\psi_1)}\biggr) = \frac{\mathbf{1}_{\langle
    d\rangle}}{\prod_{i=0}^n \prod_{\substack{b \colon \langle b
      \rangle = \langle dw_i \rangle \\ 0<b\leq dw_i}}(w_iP + bz)}.
\]
\end{proof}

Note that $M_d$ can consist of several connected components, some of
which can be singular or of excess dimension; this does not affect the
calculation.  Similarly the graph space $G_d$ also, in general, has
several irreducible or connected components.  The fact that $u$ is
virtually birational implies that only the component which generically
consists of morphisms from $\p^{1,r}$ contributes nontrivially to the
calculation.

\subsection{The Stack $L_d$ of Polynomial Maps and the Morphism $j
  \colon \p(V^f) \hookrightarrow L_d$}
\label{4sec:space-polyn-maps}

\begin{dfn}
  Let $\CC(w)$ denote the one-dimensional vector space $\CC$ equipped
  with a weight-$w$ action of $\Cstar$.
\end{dfn}

By $\Cstar$ here we mean the $\Cstar$ which occurs in the quotient
\eqref{2eq:defofwps} not the $\Cstar$ which acts by ``rotation of
loops'' \eqref{4eq:action}.  Recall that $d=m/r$ in lowest terms, and
that $f = \langle -d \rangle \in F$.

\begin{dfn}
  \label{4dfn:1}
  \[
  L_d =\Biggl[ \biggl( 
  \bigoplus_{i=0}^n  \C(-w_i)^{\oplus \bigl( 1+\lfloor dw_i \rfloor
    \bigr)}  
  - \{ 0 \} 
  \biggr)   / \Cstar \Biggr]
  \]
\end{dfn}

We regard $L_d$ as the stack of polynomial maps $\p^{1,r}\to
\p^\mathbf{w}$ of degree $d$, as follows.  Such a map is
given\footnote{See also Claim~\ref{4cla:L_univ} below.} by
polynomials
\[
P_0,P_1,\dots,P_n,
\] 
not all zero, where $P_i=P_i(s_0,s_1)$ is homogeneous of degree $mw_i$
in the variables $s_0,s_1$, with $\deg s_0=1, \deg s_1=r$.  Each $P_i$
can be written as
\begin{equation}
  \label{4eq:1}
P_i(s_0,s_1)= A_0s_0^{mw_i}+A_1s_0^{mw_i-r}s_1 + \cdots + 
A_{\lfloor dw_i\rfloor }s_0^{r\langle dw_i\rangle}
s_1^{\lfloor dw_i\rfloor}  
\end{equation}
and hence
\[
L_d =\Biggl[ \biggl( 
\bigoplus_{i=0}^n  \C(-w_i)^{\oplus \bigl( 1+\lfloor dw_i \rfloor
  \bigr)}  
- \{ 0 \} 
\biggr)   / \Cstar \Biggr]
\]
The stack $L_d$ is itself a weighted projective space.  Recall that
\[
V^f = \bigoplus_{i \colon fw_i\in \z} \C(-w_i)
\]
and note that $fw_i$ is an integer if and only if $dw_i$ is an
integer.  We define the map $j\colon \p (V^f) \hookrightarrow L_d$ by
\[
j\colon \C(-w_i) \ni A_i\mathbf{e}_i \mapsto A_i s_1^{dw_i}\in \C(-w_i).
\] 
The action \eqref{4eq:action} of $\C^\times$ on $\PP^{1,r}$ induces an
action on $L_d$ in the obvious way.

\begin{rem}
  \label{4rem:1}
  It is clear that $j\colon \p(V^f)\hookrightarrow L_d$ is a component
  of the $\Cstar$-fixed substack.
\end{rem}

\begin{rem} \label{4rem:fixedsubstack}
  Consider an action $\Psi \colon G\times \mathcal{X} \to \mathcal{X}$
  of a group scheme $G$ on a stack $\mathcal{X}$. A substack $\iota
  \colon \mathcal{Y}\hookrightarrow \mathcal{X}$ is \emph{fixed} by the action
  if for all schemes $S$ we have a diagram:
  \begin{equation}
    \label{4eq:4}
    \xymatrix{
      {G}(S)\times
      {\mathcal{Y}}(S)\ar[rr]^{{\text{pr}}_2(S)}
      \ar[d]_{{\text{id}}_G(S)\times
        {\iota}(S)} & &
      {\mathcal{Y}}(S)\ar[d]^{{\iota}(S)}\\ 
      {G}(S)\times {\mathcal{X}}(S)\ar[rr]_{{\Psi}(S)}
      \ar@{=>}[rru] & &
      {\mathcal{X}}(S)
    }  
  \end{equation}
  where the $\Rightarrow$ means that there is an \emph{isomorphism of
    functors}
  \begin{equation}
    \label{4eq:5}
    {\Psi}(S) \circ \bigl({\text{id}}_G(S)\times
    {\iota}(S)\bigr)
    \Rightarrow
    {\iota}(S) \circ {\text{pr}}_2(S)  
  \end{equation}
  By definition, a fixed substack $\iota \colon
  \mathcal{Y}\hookrightarrow \mathcal{X}$ is \emph{the $G$-fixed
    substack} if it satisfies the obvious universal property: if
  $j\colon \mathcal{Z} \hookrightarrow \mathcal {X}$ is any other
  fixed substack, than it factors uniquely through $\iota \colon
  \mathcal{Y} \hookrightarrow \mathcal{X}$.
\end{rem}

\begin{lem}
  \label{4lem:1}
  Let $N_j$ be the normal bundle of the inclusion $j\colon
  \p(V^f)\hookrightarrow L_d$.  The $\Cstar$-equivariant Euler class
  of $N_j$ is
  \[
  e(N_j)=\prod_{i=0}^n \prod_{\substack{b \colon \langle b \rangle =
      \langle dw_i \rangle \\0<b\leq dw_i}}(w_iP+bz).
  \]
\end{lem}
\begin{proof}
  Contemplate the following diagram on $\p (V^f)$.  The bottom two
  rows are the Euler sequence for weighted projective space:
  \[
  \xymatrix{
    & & 0 & 0 & \\
    & & \bigoplus_{i=0}^n \bigoplus_{\substack{b \colon \langle b
        \rangle = \langle dw_i \rangle \\0<b\leq dw_i}}
    \o(w_iP+bz)\ar@{=}[r]\ar[u] &
    N_j \ar[u] & \\
    0\ar[r] & \C \ar[r] &\bigoplus_{i=0}^n \bigoplus_{\substack{b \colon
        \langle b \rangle = \langle dw_i \rangle \\0\leq b\leq dw_i}}
    \o(w_iP+bz)
    \ar[r]\ar[u]& T_{L_d}\big|_{\p (V^f)} \ar[r]\ar[u]& 0\\
    0\ar[r] & \C \ar[r] \ar@{=}[u]& \bigoplus_{i:dw_i\in \z} \o(w_iP)
    \ar[r]\ar[u]& T_{\p (V^f)} \ar[r]\ar[u]& 0\\
    & & 0 \ar[u] & 0\ar[u] & }
  \]
\end{proof}

\subsection{Deformations and Obstructions}
\label{4sec:deform-stable-maps}

We review the canonical obstruction theories on $M_d$ and $G_d$ and
prove that the obstruction theory on $M_d$ is inherited from the
obstruction theory on $G_d$.

\begin{nt}
  Given a stack $\cX$ and a scheme $S$, we write $\cX(S)$ for the
  category of morphisms from $S$ to $\cX$.
\end{nt}

\subsubsection{The $\Cstar$-Action on $G_d$}
\label{4sec:cstar-action-g_d}

The (left) action of $\Cstar$ on $\p^\mathbf{w}\times \p^{1,r}$, where
$\Cstar$ acts on the second factor only via \eqref{4eq:action},
induces an action on the stack $G_d$ by ``dragging'' the image of the
morphism.  More precisely, given a scheme $S$, an object of
${G}_d (S)$ is a stable morphism over $S$:
\[
\xymatrix{
\mathcal{C} \ar[rr]^f \ar[d]^p & & \p^\mathbf{w}\times\p^{1,r} \\
S\ar@/^1pc/[u]^{\sigma} & &\\}
\]
and the group action is described as
\begin{align*}
\lambda \colon f \mapsto {}^\lambda f = \ell_{\lambda^{-1}} \circ f  &&
\lambda \in \Cstar
\end{align*}
where $\ell_\lambda \colon \p^\mathbf{w}\times\p^{1,r} \to \p^\mathbf{w}\times\p^{1,r}$ is left
translation by $\lambda$. 

\subsubsection{The Stack $\iota \colon M_d \hookrightarrow G_d$ is
  Part of the Fixed Substack}
\label{4sec:iota-colon-m_d}

We now construct the morphism of stacks $\iota \colon M_d \to G_d$
used in Theorem~\ref{4thm:1}.  Throughout this subsection we assume
that $\langle d \rangle \ne  0$.  The results remain true if $\langle
d \rangle = 0$; the proofs in this case are slightly different but
similar and easier.

For all schemes $S$, we
need functors ${\iota}(S)\colon {M}_d(S)\to
{G}_d(S)$ satisfying various compatibilities. An object of
${M}_d(S)$ is a stable morphism:
\begin{equation}
  \label{4eq:objectofMdS}
  \xymatrix{\mathcal{C}^\prime\ar[d]^{p^\prime} 
    \ar[rr]^{f^\prime} & & \p^\mathbf{w}\\
    S\ar@/^1pc/[u]^{\sigma^\prime} & & 
  }
\end{equation}
where $\sigma^{\prime}$ is a section of the gerbe marking. Denote by
$C_{r,r}$ the twisted curve with coarse moduli space $\p^1$ and stack
structure with stabilizer $\mu_r$ at $0, \infty$ determined by
charts\footnote{We have $C_{r,r} = \left[\PP^1/\mu_r\right]$
  where $\mu_r$ acts via $\xi\colon \thinspace[a_0:a_1] \mapsto [\xi
  a_0: a_1]$.}
\begin{align*}
&  [\C /\mu_r] \quad \text{where}\;\mu_r\; \text{acts in the standard way
  at}\; 0, \; \text{and} \\
&  [\C /\mu_r] \quad \text{where}\;\mu_r\; \text{acts as}\; \zeta \colon z
\mapsto \zeta^{-1} z \; \text{at}\; \infty.
\end{align*}
There is a natural
morphism of stacks $C_{r,r}\to \p^{1,r}$ of degree $1/r$; this morphism is
representable at $0$ and nonrepresentable at $\infty$. We denote by
\[
\xymatrix{\mathcal{C}^{\prime\prime} \ar[d]^{p^{\prime\prime}}
  \ar[rr]^{f^{\prime\prime}} & & \p^{1,r}\\ 
S\ar@/^1pc/[u]^{\sigma^{\prime\prime}_0,\; \sigma^{\prime \prime}_\infty} & & 
}
\] 
the trivial family $\mathcal{C}^{\prime\prime}=S\times C_{r,r}$ over
$S$ with (nonrepresentable) morphism to $\p^{1,r}$. By definition, the
functor ${\iota} (S) \colon {M}_d (S)\to {G}_d (S)$ maps the family
\eqref{4eq:objectofMdS} to the family
\begin{equation}
  \label{4eq:3}
\mathcal{C}/S\to \p^\mathbf{w}\times\p^{1,r} :=
\xymatrix{\mathcal{C}^\prime\cup_{\sigma^\prime,
    \sigma^{\prime\prime}_\infty} \mathcal{C}^{\prime
    \prime}\ar[d]^{p^\prime \cup p^{\prime\prime}}
  \ar[rrr]^{(f^\prime,\infty)\cup (ev_1^\prime
    p^{\prime\prime},f^{\prime \prime})} & & &\p^\mathbf{w}\times\p^{1,r}\\
  S\ar@/^1pc/[u]^{\sigma^{\prime\prime}_0} & & &}  
\end{equation}
The glued family $\mathcal{C}^\prime\cup_{\sigma^\prime,
  \sigma^{\prime\prime}_\infty} \mathcal{C}^{\prime \prime}$ here is
constructed using \cite{AGV:2}*{Proposition~A.0.2}. It is easy to see
that the functors ${\iota}(S) \colon {M}_d (S)\to {G}_d(S)$ combine to
give a closed substack $\iota \colon M_d\hookrightarrow G_d$.

\begin{lem}
  \label{4lem:2} The substack
  $\iota \colon M_d \hookrightarrow G_d$ is a $\Cstar$-fixed substack.
\end{lem}

\begin{proof}[Sketch of proof] This is an extended exercise in
  unravelling the definition of fixed substack, which was given in
  Remark~\ref{4rem:fixedsubstack}.  We give a sketch since we could
  find no adequate reference in the literature.  A well-written and
  careful treatment of group actions on stacks can be found in
  \cite{Romagny}.

  Consider an object $\xi_S^\prime = (f^\prime \colon
  \mathcal{C}^\prime/S \to \p^\mathbf{w})$ of ${M}_d(S)$ as in
  \eqref{4eq:objectofMdS} and let
  \[
  \xi_S={\iota}(S)(\xi^\prime_S)=(f\colon \mathcal{C}/S \to
  \p^\mathbf{w}\times \p^{1,r})
  \]
  be the family of diagram~(\ref{4eq:3}).  We must exhibit, for every
  $S$-point $\lambda \in \Mor (S, \Cstar)$, an isomorphism from
  ${}^\lambda \xi_S$ to $\xi_S$ which is sufficiently natural that it
  satisfies all the necessary compatibilities and produces the
  isomorphism of functors $\Rightarrow$.  This all follows from:

  \begin{cla}
    \label{4cla:equiv}
    In the notation of the preceding paragraph, there is a natural
    $\Cstar$-action on $\mathcal{C}$ which covers the trivial action on
    $S$ such that the morphism $f\colon \mathcal{C}\to
    \p^\mathbf{w}\times \p^{1,r}$ is $\Cstar$-equivariant.
  \end{cla}
  
  This is obvious: the family $\mathcal{C}$ is obtained by glueing the
  families $\mathcal{C}^\prime$ and
  $\mathcal{C}^{\prime\prime}=S\times C_{r,r}$.  $\Cstar$ acts on
  $\mathcal{C}^{\prime \prime}$ by acting on the second factor alone,
  and this action glues with the trivial action on
  $\mathcal{C}^\prime$ to give an action on $\mathcal{C}$.

  Now the Claim precisely says that, for all $\lambda \in \Cstar
  (S)$, the left translation $\ell_{\lambda^{-1}} \colon
  \mathcal{C}\to \mathcal{C}$ sits in a commutative diagram:
  \[
  \xymatrix{\mathcal{C}\ar[dr]_{{^\lambda} f} \ar[rr]^{\ell_{\lambda^{-1}}} & &
    \mathcal{C}\ar[dl]^f\\
    & \p^\mathbf{w} \times \p^{1,r} }
  \]
  That is, exactly as desired, $\ell_{\lambda^{-1}}$ defines an isomorphism
  from ${}^\lambda \xi_S$ to $\xi_S$.  This shows that $\iota \colon M_d
  \hookrightarrow G_d$ is a $\Cstar$-fixed substack.
\end{proof}

We show in Lemma~\ref{4lem:9} below that $\iota \colon M_d \hookrightarrow G_d$ is a
\emph{part} of the $\Cstar$-fixed substack.

\subsubsection{Perfect Obstruction Theory}
\label{4sec:obstruction-theory}

We recall some facts about perfect obstruction theories from
\cite{Behrend--Fantechi, Li--Tian}. For a morphism $q\colon
\mathcal{X} \to \mathcal{S}$ of stacks we denote by $L_q^\bullet$ the
first-two-term cutoff of the cotangent complex of $q$.  
The official reference for the cotangent complex is
\cite{Illusie:1,Illusie:2}, but an accessible introduction to the
first-two-term cutoff can be found in \cite{Grothendieck}.  Recall that a
relative perfect obstruction theory is a $q$-perfect 2-term complex
$E^\bullet$ on $\mathcal{X}$ together with a morphism $\varphi \colon
E^\bullet \to L_q^\bullet$ which is an isomorphism on $H^0$ and
surjective on $H^{-1}$; a relative perfect obstruction theory produces
a virtual fundamental class $\mathbf{1}^\text{vir}_q \in CH_\bullet
\(\mathcal {X}\)$.
    
Let $\mathcal{X}$ be a stack and $d\in H_2(\mathcal{X};\QQ)$.  Denote,
as usual, by $\mathcal{X}_{g,n,d}$ the moduli stack of genus-zero
$n$-pointed stable morphisms to $\cX$ of degree $d$; analogous remarks
apply to the stacks $\Mbar_{g,n}(\mathcal{X},d)$ of $n$-pointed stable
morphisms with sections to all gerbes.  There are, as we now recall,
two natural obstruction theories on $\mathcal{X}_{g,n,d}$ and they
produce the same virtual fundamental class.  Consider the universal
family:
\[
\xymatrix{
  \mathcal{U} \ar[r]^f \ar[d]_\pi & \mathcal{X} \\
  \mathcal{X}_{g,n,d} & }
\]
\begin{enumerate}
\item The \emph{relative obstruction theory} $E^{\bullet \,
    \vee}_{\text{rel}}=R\pi_\star f^\star T_\mathcal{X}$ is an
  obstruction theory relative to the canonical morphism $q\colon
  \mathcal{X}_{g,n,d} \to \mathfrak{M}^\text{tw}_{g,n}$ to the Artin
  stack of pre-stable twisted curves. The relative obstruction theory
  is used in \cite{AGV:1, AGV:2}, because it is well-suited to
  checking the axioms of Gromov--Witten theory.
\item The \emph{absolute obstruction theory} is
\[
E^{\bullet \, \vee} = R\pi_\star R\underline{\Hom}_{\o_\mathcal{U}} \bigl(
L_f^\bullet,\o_{\mathcal{U}} \bigr),
\quad
\text{where}
\quad
L_f^\bullet =[f^\star\Omega^1_\mathcal{X} \to \Omega^1_\pi (\log)]
\]
is the cotangent complex of $f$; here $\Omega^1_\pi (\log)$ denotes
the sheaf of K\"ahler differentials with logarithmic poles along the
markings.
\end{enumerate}
It is well-known that the absolute and relative obstruction theories
produce the same fundamental class (see \citelist{
  \cite{Graber--Pandharipande}*{Appendix~B}
  \cite{Okounkov--Pandharipande}*{Proposition~5.3.5}
  \cite{Kim--Kresch--Pantev}*{Proposition~3} 
}). In what follows, we
use the absolute theory.

\subsubsection{Obstructions and Virtual Normal Bundle}
\label{4sec:obstr-virt-norm}

In this section we compare obstruction theories and calculate the
virtual normal bundle of $\iota \colon M_d \hookrightarrow G_d$.

We recall a few general notions from \cite{Graber--Pandharipande}.
Let $G$ be a group scheme acting on a stack $\mathcal{X}$ and let
$E^\bullet\to L^\bullet$ be a $G$-linearized perfect obstruction
theory. Let $\iota \colon \mathcal{Y}\hookrightarrow \mathcal{X}$ be the
$G$-fixed substack. Then $G$ acts on
$\left.E^{\bullet}\right|_{\mathcal{Y}}$, and it is a fact that the
complex of $G$-invariants $\left.E^{-1}\right|_{\mathcal{Y}}^G \to
\left.E^{0}\right|_{\mathcal{Y}}^G$ is an obstruction theory for
$\mathcal{Y}$.  We call this the \emph{inherited obstruction theory}.
Writing $\left.E^{i}\right|_{\mathcal{Y}} =
\left.E^{i}\right|_{\mathcal{Y}}^G +
\left.E^{i}\right|_{\mathcal{Y}}^\text{mov}$, the moving part
$\left.E^{0}\right|_{\mathcal{Y}}^{\text{mov}\,\vee} \to
\left.E^{-1}\right|_{\mathcal{Y}}^{\text{mov}\,\vee}$ is the
\emph{virtual normal bundle}.

\begin{lem}
  \label{4lem:3} \ 
  \begin{enumerate}
  \item The obstruction theory on $M_d$ inherited from $\iota \colon
    M_d \hookrightarrow G_d$ is the natural absolute obstruction theory on
    $M_d$.
  \item Denoting by $N^\text{vir}_\iota$ the virtual normal bundle of
  $\iota$, we have 
  \[
  e(N_\iota^\text{vir})= z(z-\psi_1).
  \]
  \end{enumerate}
\end{lem}
 
\begin{proof}[Sketch of proof] The statement is well-known in a
  similar context, so we just give a sketch of the proof here. We
  start with an object $f^\prime\colon \mathcal{C}^\prime/S \to
  \p^\mathbf{w}$ of ${M}_d (S)$ as in \eqref{4eq:objectofMdS} and apply
  the functor ${\iota}(S)$ to make $f\colon \mathcal{C}/S \to
  \p^\mathbf{w}\times \p^{1,r}$ as in diagram~(\ref{4eq:3}). The first
  statement means that the natural homomorphism
  \begin{equation}
    \label{4eq:8}
    Rp_\star
    R\underline{\Hom}_{\o_\mathcal{C}} (L_f^\bullet, \o_\mathcal{C}) 
    \to
    Rp^\prime_\star R\underline{\Hom}_{\o_{\mathcal{C}^\prime}}
    (L_{f^\prime}^\bullet, \o_{\mathcal{C}^\prime}) 
  \end{equation}
  induces an isomorphism from the direct summand $Rp_\star^{\Cstar}
  R\underline{\Hom}_{\o_\mathcal{C}} (L_f^\bullet, \o_\mathcal{C})$ to
  $ Rp^\prime_\star R\underline{\Hom}_{\o_{\mathcal{C}^\prime}}
  (L_{f^\prime}^\bullet, \o_{\mathcal{C}^\prime})$. Since both
  complexes are \emph{perfect}, we can check this after base change to
  all geometric points; in effect we can and do from now on assume
  that $S=\Spec \C$, that $\mathcal{C}=C$ is a pre-stable curve over
  $\Spec \C$, etc.

Applying the cohomological functor $R\Hom_{\o_C}(-, \o_C)$ to the exact
triangle
\[
\Omega^1_p(\log) \to L_f^\bullet \to f^\star \Omega^1_{\p^\mathbf{w}\times
  \p^{1,r}}[1]\overset{+1}{\longrightarrow}
\]
we calculate $E^{0\,\vee}=\mathbb{T}^1_f$ and
$E^{-1\,\vee}=\mathbb{T}^2_f$ from the well-known exact
sequence:
\begin{multline}
  \label{4xact:defo}
  0 \to H^0(C, \Theta_C(-\log)) \to H^0(C, f^\star T_{\p^\mathbf{w}\times
    \p^{1,r}}) \to \mathbb{T}^1_f \to \\
  \to \Ext^1_{\o_C}(\Omega^1_C(\log), \o_C) \to H^1(C, f^\star
  T_{\p^\mathbf{w}\times \p^{1,r}})\to \mathbb{T}^2_f \to 0.
\end{multline}
Our goal is to determine each piece in the exact
sequence~(\ref{4xact:defo}) as a representation of $\Cstar$; we make
the following simple observations:
\begin{enumerate}
\item $\Theta_C(-\log)=\Theta_{C^\prime}(-\log)\oplus \Theta_{\p^1}(-0-\infty)$, hence
  \[
  H^0(C, \Theta_C(-\log))=H^0(C^\prime,
  \Theta_{C^\prime}(-\log))\oplus \C(z)
\]
with the first summand a trivial representation.
\item $f^\star T_{\p^\mathbf{w}\times \p^{1,r}} = f_1^\star
  T_{\p^\mathbf{w}} \oplus f_2^\star T_{\p^{1,r}}$, where $f_1\colon C
  \to \p^\mathbf{w}$ and $f_2\colon C \to \p^{1,r}$ are the natural
  morphisms.  Thus
\[
H^0(C, f^\star T_{\p^\mathbf{w}\times \p^{1,r}})=
H^0(C^\prime, f^{\prime\,\star} T_{\p^\mathbf{w}})\oplus H^0(\p^{1,r},
T_{\p^{1,r}}) 
\]
where the first summand is $\Cstar$-fixed and the second summand is
moving (and easy to calculate as a representation using the
equivariant Euler sequence on $\p^{1,r}$).
\item We calculate $\Ext^1_{\o_C}(\Omega^1_C(\log), \o_C)$ with the
  standard local-to-global spectral sequence:
  \begin{multline*}
    0 \to H^1(C,\Theta_C(-\log))\to \Ext^1_{\o_C}(\Omega^1_C(\log), \o_C)\to
    \\ \to H^0\bigl(C, \underline{\Ext}^1_{\o_C}(\Omega^1_C(\log),
    \o_C)\bigr) \to 0
  \end{multline*}
Now $H^1(C,\Theta_C(-\log))= H^1(C^\prime, \Theta^\prime(-\log))$ is a
trivial representation, whereas 
\begin{multline*}
  H^0\bigl(C, \underline{\Ext}^1_{\o_C}(\Omega^1_C(\log), \o_C)\bigr)=\\
  =H^0\bigl(C^\prime,
  \underline{\Ext}^1_{\o_{C^\prime}}(\Omega^1_{C^\prime}(\log),
  \o_{C^\prime})\bigr)\oplus \bigl(T_{C^\prime, \sigma^\prime}\otimes
  \C(z)\bigr)
\end{multline*}
where the first summand is the trivial representation and the second
summand is isomorphic to $\C(z)$. From this and the 5-lemma, we
conclude that
\[
\Ext^1_{\o_C}(\Omega^1_C(\log), \o_C)=
\Ext^1_{\o_{C^\prime}}(\Omega^1_{C^\prime}(\log), \o_{C^\prime})
\oplus \bigl(T_{C^\prime, \sigma^\prime}\otimes \C(z)\bigr)
\]
as the sum of fixed and moving parts.
\item As before, 
\[
H^1(C, f^\star T_{\p^\mathbf{w}\times \p^{1,r}})=
H^1(C^\prime, f^{\prime\,\star} T_{\p^\mathbf{w}}).
\]
\end{enumerate}
Using the above facts and the 5-lemma it is easy to finish the proof.
\end{proof}

\subsection{Construction and Properties of the Morphism $u$}
\label{4sec:construction-map-u}

We give a precise construction of the morphism $u$ following closely
the argument of Jun Li \cite[Lemma~2.6]{Lian--Liu--Yau}. Finally, we show
that the morphism $u\colon G_d \to L_d$ is virtually birational.

\begin{lem}
  \label{4lem:4}
  There is a natural morphism $u\colon G_d \to L_d$.
\end{lem}

\begin{proof}
  We sketch the proof, which follows closely
  \cite[Lemma~2.6]{Lian--Liu--Yau}.  For all schemes $S$, we construct
  functors ${G}_d(S)\to {L}_d(S)$. This is not
  difficult to do since $L_d$ is itself a weighted projective space.
  It therefore satisfies a universal property which makes it easy to
  construct elements of ${L}_d(S)$. Let us spell this out
  more precisely. We denote:
  \begin{align*}
    W=\C(-1) \oplus \C(-r), 
    &&
    \text{so that}
    &&
    \p^{1,r}=\left[ W - \{0\} / \mathbb{T}^1 \right].
  \end{align*}
  Note that the free polynomial algebra $S^\ast W^\vee$ generated by
  $W^\vee$ is a representation of $\Cstar$.  We denote by $S^m W^\vee$
  the isotypic component on which $\Cstar$ acts with weight $m\in \z$;
  $S^\ast W^\vee$ is generated by a basis element $s_0\in W^\vee$ of
  degree $1$ and a basis element $s_1\in W^\vee \cap S^r W^\vee$ of
  degree $r$.  A polynomial map $\p^{1,r}\to \p^\mathbf{w}$ of degree
  $d=m/r$ is given by polynomials $P_0,\dots, P_n \in S^{mw_i}
  W^\vee$, not all identically zero:
  \[
  L_d = \biggl[ \Bigl( \oplus_{i=0}^n S^{mw_i}W^\vee - \{0\} \Bigr)/\Cstar \biggr]. 
  \]
From this we conclude:
\begin{cla}
  \label{4cla:L_univ}
  Let $S$ be a scheme.  An object of ${L}_d(S)$ consists of a
  line bundle $\mathcal{L}$ on $S$ and a nowhere vanishing sheaf
  homomorphism:
\[
(P_0,\dots,P_n)\colon \o_S^n \to \oplus_{i=0}^n \mathcal{L}^{\otimes
  w_i} \otimes S^{mw_i} W^\vee .
\]
\end{cla}
Let us now proceed to the proof of Lemma~\ref{4lem:4}. 
An object of ${G}_d (S)$ is a stable morphism:
\begin{equation}
  \label{4eq:6}
\xymatrix{
\mathcal{C} \ar[rr]^{(p_2,p_3)} \ar[d]_{p_1} & & \p^\mathbf{w}\times
\p^{1,r} \\
S & & }  
\end{equation}
(Depending on whether or not $d$ is an integer, there may be a
section $\sigma \colon S \to \mathcal{C}$; the section plays no role
in what follows.) Let us rearrange the diagram as:
\begin{equation}
  \label{4eq:7}
\xymatrix{
  \mathcal{C} \ar[rr]^{p_2} \ar[d]_{q:=(p_1,p_3)} & & \p^\mathbf{w} \\
  S \times \p^{1,r} & &}
\end{equation}
\begin{cla} \ 
  \begin{enumerate}
  \item The sheaves $\mathcal{F}_k=q_\star p_2^\star
    \o_{\p^\mathbf{w}}(k)$ are flat over $S$ and generically of rank 1.
  \item There is a line bundle $\mathcal{L}$ on $S$ such that
    \[
    \det \mathcal{F}_k = \mathcal{L}^{\otimes k} \boxtimes 
    \o_{\p^{1,r}}(mk).
    \]
  \end{enumerate}
\end{cla}
This is proved in \cite{Lian--Liu--Yau}, and it easily implies the
result.  The canonical sections $x_i\in H^0(\p^\mathbf{w}, \o(w_i))$
give elements $p_2^\star x_i \in H^0(S \times \p^{1,r},
\mathcal{F}_{w_i})$, and using the canonical sheaf homomorphism
$\mathcal{F}_k \to \det \mathcal{F}_k$ ($\mathcal{F}_k$ has rank 1!),
these map to elements $P_i$ of
\[
  H^0(S \times \p^{1,r},\mathcal{L}^{\otimes w_i} \boxtimes
  \o_{\p^{1,r}}(mw_i))
  =H^0(S, \mathcal{L}^{\otimes w_i}\otimes S^{mw_i} W^\vee).
\]
Thus we have constructed a sequence $(P_0,\ldots,P_n)$ of elements of
$H^0(S, \mathcal{L}^{\otimes w_i}\otimes S^{mw_i} W^\vee)$ and this,
by virtue of Claim~\ref{4cla:L_univ}, gives an object of ${L}_d(S)$.
\end{proof}

It is useful to know the morphism $u$ explicitly at geometric points.
Consider an element $\varphi \colon \cC \to \p^\mathbf{w}\times
\p^{1,r}$ of $\Mbar_{0,1} (\p^\mathbf{w}\times \p^{1,r}, d\times
\frac1{r})$. Write
\begin{equation}
  \label{4eq:curvesplit}
  \cC=\bigcup_{j=0}^N \cC_j
\end{equation}
where
 \begin{enumerate}
 \item $\cC_0$ is the distinguished component mapping 1-to-1 to
   $\p^{1,r}$; and
 \item the curves $\cC_j$ for $j\geq 1$ are ``vertical'': they map to
   points $y_j \in \PP^{1,r}$ given by equations $s_1-a_j s_0^r=0$.
 \end{enumerate}
 \emph{Assume for simplicity that the marked point $x_0\in \cC$ lies
   on $\cC_0$}, and note that:
 \begin{enumerate}
 \item the marked point $x_0$ lies above $0\in \p^{1,r}$;
 \item for each $j\geq 1$ the curve $\cC_j$ meets $\cC_0$ in a unique
   point $x_j$, which lies above $y_j \in \p^{1,r}$, and the induced
   morphism $\rho_j \colon (\cC_j,x_j) \to \p^\mathbf{w}$ is
   representable and stable;
 \item the morphism $\rho_0 \colon \bigl(\cC_0,
   \{x_0,x_1,\ldots,x_N\}\bigr) \to \p^\mathbf{w}$ is representable
   and pre-stable.
\end{enumerate}
Write $d_j =\deg \rho_j$ and $\overline{f_j} = \langle -d_j\rangle$, so
that $\rho_j \in \Mbar_{0,1}(\p^\mathbf{w}, d_j)$.  It is clear that
$d=\sum_{j=0}^N d_j$.  The gerbe at $x_j$ evaluates to $\p (V^{f_j})$
where $f_0=f=\langle -d \rangle$ and, for $j\geq 1$,
\[
f_j = 
\begin{cases}
  1 - \overline{f_j} &\text{if $\overline{f_j} \not = 0$} \\
  0                      &\text{if $\overline{f_j}=0$.}
\end{cases}
\]
\begin{lem}
  \label{4lem:8}
  In these circumstances, the polynomial map $u(\varphi) \in L_d$
  constructed in Lemma~\ref{4lem:4} is given by homogeneous
  polynomials:
  \[
  \begin{pmatrix}
    P_0(s_0,s_1) \\
    \vdots      \\
    P_i(s_0,s_1) \\
    \vdots      \\
    P_n(s_0,s_1)
  \end{pmatrix}
  =
  \begin{pmatrix}
    Q_0(s_0,s_1) \prod_{j=1}^{N} (s_1-a_j s_0^r)^{\lfloor
      d_j\rfloor w_0}\\
    \vdots      \\
    Q_i(s_0,s_1) \prod_{j=1}^{N} (s_1-a_j s_0^r)^{\lfloor
      d_j\rfloor w_i}\\
    \vdots      \\
    Q_n(s_0,s_1) \prod_{j=1}^{N} (s_1-a_j s_0^r)^{\lfloor d_j\rfloor
      w_n}
  \end{pmatrix}
\]
where $\deg Q_i=r\bigl(d_0+\sum_{j=1}^N f_j\bigr)w_i$.  
\end{lem}

\begin{proof}
  This follows closely the classical case
  \cite[Lemma~2.6]{Lian--Liu--Yau}.
\end{proof}

We have
\begin{align*}
\deg P_i& =r\Bigl(d_0+\sum_{j=1}^N f_j\Bigr)
w_i+r\Bigl(\sum_{j=1}^N \lfloor d_j\rfloor\Bigr)w_i\\
& =r\Bigl(d_0+\sum_{j=1}^N (\lfloor d_j\rfloor + f_j) 
\Bigr)w_i  \\ &=rdw_i  =mw_i.  
\end{align*}
In addition, one should note that the polynomials $Q_i$ themselves
usually must contain common factors which account for the ``stacky
behaviour'' of the morphism $\rho_0$ above the points $y_j \in
\p^{1,r}$. More precisely, for all $i$:
\[
\text{$(s_1-a_j s_0^r)^{\langle \overline{f_j}w_i\rangle + f_jw_i}$
 is a factor of $Q_i(s_0,s_1)$}
\] 
and it is an exact factor for at least one $i$ such that $w_if_j$ is
an integer.

\begin{cor}
  \label{4cor:2}
  The basic diagram of Theorem~\ref{4thm:1}, where all stacks and
  morphisms have by now been constructed, is a commutative diagram of
  stacks with $\Cstar$-action. \qed
\end{cor}

\begin{lem}
  \label{4lem:9}
  The substack $\iota \colon M_d \hookrightarrow G_d$ is the part of
  the $\Cstar$-fixed substack that lies above $j\colon
  \p(V^f)\hookrightarrow L_d$.
\end{lem}

\begin{proof} 
  The basic diagram of Theorem~\ref{4thm:1} is a commutative diagram
  of stacks with $\Cstar$-action.  The $\Cstar$-fixed substack of
  $G_d$ is therefore a disjoint union of parts lying above the
  connected components of the $\Cstar$-fixed substack of $L_d$.
  $j\colon \p(V^f)\hookrightarrow L_d$ is one of these components, and we show
  that $\iota \colon M_d\hookrightarrow G_d$ is the part of the $\Cstar$-fixed
  stack lying above $\p(V^f)$ by showing that it has the required
  universal property.

  First, we show that this is so over geometric points. Let $\varphi
  \colon \cC \to \p^\mathbf{w}\times \p^{1,r}$ be a $\Cstar$-fixed
  point of $G_d$. Write $\cC=\bigcup_{j=0}^N\cC_j$ as in
  \eqref{4eq:curvesplit}, so that $\cC_0$ is the distinguished
  component mapping 1-to-1 to $\p^{1,r}$ and the $\cC_j$ are vertical
  for $j\geq 1$. Since $\varphi$ is $\Cstar$-fixed, by the very way
  the $\Cstar$-action is defined, the image $\varphi (\cC)\subset
  \p^\mathbf{w}\times \p^{1,r}$ is invariant under the action of
  $\Cstar$ on $\p^\mathbf{w} \times \p^{1,r}$ acting on the second
  factor only.  This implies that $\varphi (\cC_0)$ is a
  \emph{horizontal curve}; it then follows from Lemma~\ref{4lem:8} and
  Corollary~\ref{4cor:2} that there is only one vertical curve $\cC_j$
  and that it is joined to $\cC_0$ over $\infty \in \p^{1,r}$.  In
  other words, $\varphi$ is \emph{isomorphic} to a point in the image
  of $\iota$.

  We are now ready to finish the proof of the Lemma.  Consider a base
  scheme $S$ and a $\Cstar$-fixed object of ${G}_d(S)$:
  \begin{equation}
    \label{4eq:10}
\xymatrix{\mathcal{C}\ar[d]^{p} 
\ar[rr]^{f} & & \p^\mathbf{w}\times
\p^{1,r} \\S\ar@/^1pc/[u]^{\sigma} & & 
}    
  \end{equation}
  All we need to show is that
  $\mathcal{C}=\mathcal{C}^\prime\cup_{\sigma^\prime,
    \sigma^{\prime\prime}_\infty} \mathcal{C}^{\prime \prime}$ as in
  diagram~(\ref{4eq:3}). First of all, by what we said on geometric
  points, family~(\ref{4eq:10}), considered as a family of pre-stable
  curves, is the pull-back from a unique morphism to the ``boundary''
  substack
  \[
  \mathfrak{M}^\text{tw}_{0,2}\times_{B\mu_r}
  \mathfrak{M}^\text{tw}_{0,1} \to \mathfrak{M}^\text{tw}_{0,1},
  \]
  where $\mathfrak{M}^\text{tw}_{g,n}$ is the \emph{smooth} Artin
  stack of pre-stable $n$-pointed twisted curves of genus $g$
  constructed in \cite{Olsson}.  That is,
  $\mathcal{C}=\mathcal{C}^\prime\cup_{\sigma^\prime,
    \sigma^{\prime\prime}_\infty} \mathcal{C}^{\prime \prime}$
  \emph{as a family of pre-stable curves}. Now
  \cite[Proposition~5.2.2]{AGV:2} implies that
  $\mathcal{C}=\mathcal{C}^\prime\cup_{\sigma^\prime,
    \sigma^{\prime\prime}_\infty} \mathcal{C}^{\prime \prime}$ as
  families of stable morphisms.
\end{proof}

\begin{lem}
  \label{4lem:5}
  The morphism $u$ is virtually birational:
  \[
  u_\star \mathbf{1}^\text{vir}_{G_d} = \mathbf{1}_{L_d}
  \]
\end{lem}

Before proving this, it is useful to calculate virtual dimension of
the two stacks:

\begin{lem}
  \label{4lem:6}
  $\dim \mathbf{1}^{vir}_{G_d}=\dim L_d = n+\sum \lfloor
  dw_i\rfloor$.
\end{lem}

\begin{proof} We calculate using the dimension formula of
  Equation~\eqref{2eq:dim}
  \begin{align*}
\dim \mathbf{1}^{vir}_{G_d}&= 1 + \dim (\p^\mathbf{w}\times \p^{1,r})
-3 -K_{\p^\mathbf{w}\times \p^{1,r}}\cdot
\Bigl(d, \textstyle \frac1{r} \Bigr) - \text{age}\\
&=1+n+1-3 + d\Bigl(\sum_{i=0}^n w_i \Bigr)
+\frac{r+1}{r} 
-\sum_{i=0}^n\langle -fw_i \rangle
-\frac{1}{r}\\
&=n+\sum_{i=0}^{n}\lfloor dw_i\rfloor.    
\end{align*}
\end{proof}

\begin{proof}[Proof of Lemma~\ref{4lem:5}]
  There is a unique component of $G_d$ generically parametrizing
  morphisms from irreducible curves and it maps generically 1-to-1 to
  $L_d$.  This component of $G_d$ is generically smooth and of the
  expected dimension; the virtual fundamental class of this component
  therefore coincides with the usual fundamental class and pushes
  forward to give the fundamental class of $L_d$.  If a component of
  $G_d$ generically parametrizes morphisms from reducible curves, it
  maps to a proper subvariety of $L_d$.
\end{proof}

\subsection{Proof of Theorem~\ref{4thm:1}}
\label{4sec:proof-theorem}

Putting together all the pieces, we have a proof of
Theorem~\ref{4thm:1}. The existence of the commutative diagram was
shown in Corollary~\ref{4cor:2}; the first statement is
Remark~\ref{4rem:1} and Lemma~\ref{4lem:1}; the second statement is
Lemma~\ref{4lem:9} and Lemma~\ref{4lem:3}; the third statement is
Lemma~\ref{4lem:5}. \qed

\section{The Small Quantum Cohomology of Weighted Projective Spaces}
\label{sec:qu_coh_wps}

In this section we prove Theorem~\ref{1thm:1}.  As was discussed in
Section~\ref{2sec:smallJ}, and as we will see rather explicitly below,
to determine the small quantum orbifold cohomology algebra of
$\stackX$ it suffices to compute the directional derivatives
\begin{align} \label{5eq:seek}
  \left.\nabla_{\phi_i} \bJ_{\stackX}(\tau)\right|_{\tau \in
    H^2(\cX;\CC) \subset \HorbX} 
  && i \in \{1,2,\ldots,N\}.
\end{align}
where $\phi_1,\ldots,\phi_N$ is a basis for $\Horb$.  We have computed
the small $J$-function $J_{\stackX}(t)$, which is the restriction of
$\bJ_{\stackX}(\tau)$ to $H^2(\stackX;\CC) \subset \Horb$:
\[
J_{\stackX}(t) = \bJ_{\stackX}(t P).
\]
This does not, \emph{a priori}, determine the directional derivatives
\[
\left.\nabla_{y} \bJ_{\stackX}(\tau)\right|_{\tau \in H^2(\cX;\CC)
  \subset \HorbX}
\]
along directions $y$ not in $H^2(\stackX;\CC)$, but it does allow us to calculate
multiple derivatives
\[
\left.\nabla_P \cdots \nabla_P \bJ_{\stackX}(\tau)\right|_{\tau = t P}
= {\partial \over \partial t} \cdots {\partial \over \partial t}
J_{\stackX}(t).
\]
We will combine these calculations with the differential equations
\eqref{2eq:smallQDEs} to determine the directional derivatives
\eqref{5eq:seek}.

Let $N = w_0 + \cdots + w_n$ and let $c_1,\ldots,c_N$ be the sequence
obtained by arranging the terms
\[
{0 \over w_0},{1 \over w_0},\ldots,{w_0-1 \over w_0},
{0 \over w_1},{1 \over w_1},\ldots,{w_1-1 \over w_1},
\ldots,
{0 \over w_n},{1 \over w_n},\ldots,{w_n-1 \over w_n}
\]
in increasing order.  Define differential operators
\[
D_j =
\begin{cases}
  \id & j=1\\
  Q^{-c_j} e^{- c_j t} \prod_{m=1}^{j-1} \(z {\partial \over \partial t} - z
  c_m\) & 1 < j \leq N.
\end{cases}
\]

\begin{lem} \label{5lem:vi}
  There exist $v_1,\ldots,v_N \in \HorbL$ such that 
  \begin{align*}
    z^{-1} D_j J_{\stackX}(t) = 
    \nabla_{v_j} \left. \bJ_{\stackX}(\tau)\right|_{\tau = t P}
    && j \in \{1,2,\ldots,N\}.
  \end{align*}
  Furthermore 
  \begin{enumerate}
  \item[(a)] $v_1 = \fun_0$;
  \item[(b)] $v_{j+1} = Q^{c_j - c_{j+1}} e^{(c_j - c_{j+1})t}  P \circ_{t P} v_j$, $1 \leq j < N$;
  \item[(c)] $v_j = \sigma_j P^{r_j} \fun_{c_j}$, $1 \leq j \leq N$, where
    \begin{align*}
      \sigma_j & = {\prod_{m \colon c_m < c_j} \( c_j - c_m \) \over
        \prod_{i=0}^n \prod_{\substack{b \colon \langle b \rangle = 
            \langle c_j w_i \rangle \\ 0 < b \leq c_j w_i}} b} \\
      \intertext{and}
      r_j & = \# \left \{ i \mid \text{$i<j$ and $c_i = c_j$} \right \}.
    \end{align*}
  \end{enumerate}
  In particular, $v_1,\ldots,v_N$ is a basis for $\Horb$.
\end{lem}

\begin{rem} \label{5rem:sequences}
  Note that the sequence $c_1,\ldots,c_N$ is
  \[
  \overbrace{f_1,\ldots,f_1}^{\dim_{f_1} + 1},
  \overbrace{f_2,\ldots,f_2}^{\dim_{f_2} + 1},
  \ldots,
  \overbrace{f_k,\ldots,f_k}^{\dim_{f_k} + 1}
  \]
  and that the sequence $\sigma_1,\ldots,\sigma_N$ is
  \[
  \overbrace{s_1,\ldots,s_1}^{\dim_{f_1} + 1},
  \overbrace{s_2,\ldots,s_2}^{\dim_{f_2} + 1},
  \ldots,
  \overbrace{s_k,\ldots,s_k}^{\dim_{f_k} + 1}
  \]
  where $f_1,\ldots,f_k$ are defined above equation
  \eqref{1eq:obviousbasis} and $s_1,\ldots,s_k$ are defined in
  Theorem~\ref{1thm:1}.
\end{rem}

\begin{proof}[Proof of Lemma~\ref{5lem:vi}]
  The string equation \cite[Theorem~8.3.1]{AGV:2} implies that
  \[
  z \nabla_{\fun_0} \bJ_{\stackX} (\tau) = \bJ_{\stackX} (\tau), 
  \]
  so we can take $v_1 = \fun_0$.  Assume that
  \[
  z^{-1} D_j J_{\stackX}(t) = 
  \nabla_{v_j} \left. \bJ_{\stackX}(\tau)\right|_{\tau = t P}.
  \]
  for some $j$ with $1 \leq j \leq N-1$.  Since
  \[
  z {\partial \over \partial t} D_j = 
  Q^{c_{j+1} - c_j} e^{(c_{j+1} - c_j)t}
  D_{j+1}
  \]
  we have
  \begin{align*}
    z^{-1} D_{j+1} J_{\stackX}(t) 
    &= Q^{c_j - c_{j+1}} e^{(c_j - c_{j+1})t} {\partial \over \partial t} 
    D_j   J_{\stackX}(t)  \\
    &= Q^{c_j - c_{j+1}} e^{(c_j - c_{j+1})t} z {\partial \over
      \partial t}
    \(  \nabla_{v_j} \left. \bJ_{\stackX}(\tau)\right|_{\tau = t P} \)
    \\
    &= Q^{c_j - c_{j+1}} e^{(c_j - c_{j+1})t}  
    \nabla_{P \circ_\tau v_j} \left. \bJ_{\stackX}(\tau)\right|_{\tau
      = t P} && \text{\emph{cf.} \eqref{2eq:smallQDEs}}.
  \end{align*}
  Thus we can take
  \[
  v_{j+1} = Q^{c_j - c_{j+1}} e^{(c_j - c_{j+1})t}  P \circ_{t P} v_j.
  \]
  By induction, this proves the existence of $v_1,\ldots,v_N$.  It
  also proves (a) and (b).

  We know that
  \[
  \nabla_{v_j} \bJ_{\stackX}(\tau) = v_j + O(z^{-1})
  \]
  and that
  \[
  \nabla_{v_j} \left. \bJ_{\stackX}(\tau)\right|_{\tau = t P}
  = {1 \over z} D_j J_{\stackX}(t),
  \]
  so to establish (c) we need to compute the coefficient of $z$ in
  \[
  D_j J_{\stackX}(t) = z \, e^{P t/z} \sum_{\substack{d \colon d \geq 0  \\
    \langle d \rangle \in F}} Q^{d - c_j} 
  e^{(d - c_j)t} \fun_{\langle d \rangle}
  {\prod_{m=1}^{j-1} \(P + \(d-c_m\)z\) \over
    \prod_{i=0}^n \prod_{\substack{b \colon \langle b \rangle = 
        \langle d w_i \rangle \\ 0 < b \leq d w_i}} 
    \( w_i P + b z \)}.
  \]
  The degree in $z$ of the denominator of the $d$th summand here is
  \[
  \lceil w_0 d \rceil +   \lceil w_1 d \rceil + \cdots + 
  \lceil w_n d \rceil,
  \]
  which is the number of fractions
  \begin{align*}
    {k \over w_i} && \text{$k \geq 0$, $0 \leq i \leq n$}
  \end{align*}
  which are strictly less than $d$.  If $d > c_j$ then this exceeds
  the degree in $z$ of the numerator and so the $d$th summand, when
  expanded as a Laurent series in $z^{-1}$, is $O(z^{-1})$.  Recall
  that $\PP(V^f)$ is a weighted projective space of dimension
  $\dim_f$.  If $d<c_j$ then, by Remark~\ref{5rem:sequences}, there
  are $\dim_d + 1$ values of $\ell$ such that $\ell \in
  \{1,2,\ldots,j\}$ and $c_\ell = d$.  This implies that the $d$th
  summand above contains a factor of
  \[
  P^{\dim_d + 1} \fun_d,
  \]
  which vanishes for dimensional reasons.  Thus only the summand
  where $d=c_j$ contributes to the coefficient of $z$:
  \[
  D_j J_{\stackX}(t) = z e^{P t/z} \fun_{ c_j }
  { P^{r_j} \prod_{m \colon c_m<c_j} \(P + \(c_j-c_m\)z\) \over
    \prod_{i=0}^n \prod_{\substack{b \colon \langle b \rangle = 
        \langle c_j w_i \rangle \\ 0 < b \leq c_j w_i}} 
    \( w_i P + b z \)} + o(z).
  \]
  The degree in $z$ of the numerator and denominator here are equal, so
  \[
  D_j J_{\stackX}(t) = z \fun_{c_j }
  {P^{r_j} \prod_{m \colon c_m<c_j} \(c_j-c_m\) \over
    \prod_{i=0}^n \prod_{\substack{b \colon \langle b \rangle = 
        \langle c_j w_i \rangle \\ 0 < b \leq c_j w_i}} 
    b } + o(z)
  \]
  and therefore $v_j = \sigma_j P^{r_j} \fun_{c_j}$, as claimed.
\end{proof}

\begin{lem} \label{5lem:vN}
  \[
  P \circ_{t P} v_{N} = {1 \over w_0^{w_0} w_1^{w_1}
    \cdots w_n^{w_n}} Q^{1-c_N} e^{(1-c_N)t} \fun_0
  \]
\end{lem}

\begin{proof}
  On the one hand
  \begin{align*}
    \left.\nabla_{P \circ_{t P} v_N} \bJ_{\stackX}(\tau)\right|_{\tau =
      tP} 
    &= z \left.\nabla_P \nabla_{v_N} \bJ_{\stackX}(\tau)\right|_{\tau =
      tP}
    && \text{\emph{cf.} \eqref{2eq:smallQDEs}} \\
    &= {\partial \over \partial t} D_N J_{\stackX}(t)\\
    \intertext{and on the other hand}
    \left.\nabla_{P \circ_{t P} v_N} \bJ_{\stackX}(\tau)\right|_{\tau =
      tP} 
    &= P \circ_{tP} v_N + O(z^{-1}),
  \end{align*}
  so we need to compute the coefficient of $z^0$ in 
  \[
    {\partial \over \partial t} D_N J_{\stackX}(t) = 
    e^{P t/z} \sum_{\substack{d \colon d \geq 0  \\
    \langle d \rangle \in F}} Q^{d - c_N} 
    e^{(d - c_N)t} \fun_{\langle d \rangle}
    {\prod_{m=1}^N \(P + \(d-c_m\)z\) \over
      \prod_{i=0}^n \prod_{\substack{b \colon \langle b \rangle = 
          \langle d w_i \rangle \\ 0 < b \leq d w_i}} 
      \( w_i P + b z \)}.
  \]
  Arguing exactly as in the proof of Lemma~\ref{5lem:vi}(c) we see
  that only the summand with $d=1$ contributes and that
  \begin{align*}
    {\partial \over \partial t} D_N J_{\stackX}(t) 
    &= Q^{1 - c_N} e^{(1 - c_N)t} \fun_0 
    { \prod_{m=1}^N \(1-c_m\) \over
      \prod_{i=0}^n w_i!} + O(z^{-1}). \\
    \intertext{Thus}    P \circ_{t P} v_{N}
    &= Q^{1 - c_N} e^{(1 - c_N)t} \fun_0 
    { \prod_{m=1}^{N} \(1-c_m\) \over
      \prod_{i=0}^n w_i!} \\
    &= {1 \over w_0^{w_0} w_1^{w_1}\cdots w_n^{w_n}}
    Q^{1 - c_N} e^{(1 - c_N)t} \fun_0.
  \end{align*}
\end{proof}

Lemma~\ref{5lem:vi} and Lemma~\ref{5lem:vN} together show that
the matrix of small orbifold quantum multiplication $P \circ_{t P}$
with respect to the the basis
\begin{equation}
  \label{eq:temporarybasis}
  Q^{c_1} e^{c_1 t} v_1, Q^{c_2} e^{c_2 t} v_2, \ldots, Q^{c_N} e^{c_N
    t} v_N
\end{equation}
is
\[
\begin{pmatrix}
  0 & 0 & 0 & \cdots & 0 & {Q e^t \over w_0^{w_0} w_1^{w_1} \cdots w_n^{w_n}}\\
1 & 0 & 0 & \cdots & 0 & 0 \\
0 & 1 & 0 & \ddots & & 0 \\
\vdots & & \ddots & & \vdots & \vdots \\
0 & & & \ddots & 0 & 0 \\
0 & 0 & \cdots & 0 & 1 & 0 \\
\end{pmatrix}.
\]

\begin{cor} \label{5cor:theorem1.1}
  Theorem~\ref{1thm:1} holds.
\end{cor}

\begin{proof}
  The basis \eqref{eq:temporarybasis} differs from the basis
  \eqref{1eq:obviousbasis} by factors of $\sigma_j$, $Q^j$, and $e^{c_j
    t}$.  Taking account of these differences yields
  Theorem~\ref{1thm:1}.
\end{proof}

\section{Weighted Projective Complete Intersections}
\label{sec:weighted_intersections}

Let $\cX$ be a quasismooth complete intersection of type $(d_0,d_1,
\dots, d_m)$ in $\p^\mathbf{w}$ and let $\iota \colon \cX \to \stackX$
be the inclusion.  Let $k_{\cX}=\sum_{j=0}^m d_j -\sum_{i=0}^nw_i$.
The main result of this section, Corollary~\ref{6cor:hyp}, determines
part of the big $J$-function of $\cX$; it applies to quasismooth
complete intersections with $k_{\cX} \leq 0$.

We begin with a combinatorial Lemma.

\begin{lem}
  \label{6lem:1} \ 
    \begin{enumerate}
    \item If $k_{\cX} \leq 0$ then for all $f\in F$,
    \[ \sum_{j=0}^m \lceil fd_j\rceil -\sum_{i=0}^n \lceil fw_i\rceil
    \leq fk_{\cX}.
    \]
  \item If $k_{\cX}= 0$ then for all non-zero $f\in F$,
    \[
    \sum_{j=0}^m \lceil fd_j\rceil -\sum_{i=0}^n \lceil fw_i\rceil
    <0.
    \]
  \end{enumerate}
\end{lem}
\begin{proof}
  The proof is elementary; see \cite[Section~8]{Fletcher} for some
  useful facts about quasismooth complete intersections. Fix $f\in F$
  and let $I=\{i \mid w_if\in \z\}$.  Since $\cX$ is quasismooth along
  $\p(V^f)\subset \p^\mathbf{w}$, we can reorder the $d_j$ and the
  $w_i$ such that:
\begin{enumerate}
\item For $j\leq l$, $fd_j$ is not an integer and there is a monomial
  $x_I^{M_I}$ in the variables $\{x_i\mid i \in I\}$ such that
  $x_jx_I^{M_I}$ has degree $d_j$; in particular, this implies that
  $fd_j \equiv fw_j \mod \z$.
\item For $l<j$, there is a monomial
  $x_I^{M_I}$ of degree $d_j$ in the variables $\{x_i\mid i \in I\}$;
  in particular, this implies that $fd_j$ is an integer.
\end{enumerate}
Then:
\begin{equation}
  \label{6eq:pain}
  \begin{split}
    \sum_{j=0}^m \lceil fd_j\rceil &= 
    fk_{\cX}+ \sum_{i=0}^l \lceil fw_i
    \rceil + \sum_{i \in I} fw_i\, + \sum_{i \not \in \{0,\dots,l\}\cup I}
    fw_i \\
    &\leq fk_{\cX} + \sum_{i=0}^n \lceil fw_i \rceil  
  \end{split}
\end{equation}
and this is part (1) of the statement. If $k_{\cX}=0$ then part (2)
also follows unless we have equality in Equation~\eqref{6eq:pain},
that is unless $\{0,\dots,l\}\cup I=\{0,\dots,n\}$. We show that this
leads to a contradiction. Let $G_0,\dots, G_m$ be the equations of
$\cX$ of degrees $\deg G_j =d_j$.  For $j=0, \dots , l$, we have that
$fd_j\not \in \z$; this implies that $\p (V^f)=\{x_0 =\cdots =
x_l=0\}$ is an irreducible component of $\{G_0=\cdots = G_l=0\}$.
This in turn implies that $\cX$ itself is reducible, a contradiction.
\end{proof}

\begin{cor} \ 
  \label{6cor:hyp}
  \begin{enumerate}
  \item If $k_{\cX} <0$, then 
    \[
    I_{\cX}(t)= \iota_\star \(z+\tau(t)+O(z^{-1})\)
    \]
    for some function $\tau \colon \C \to \HorbXL$, and 
    \[
    \iota_\star \, \bJ_{\cX} \(\tau(t)\)=I_{\cX}(t).
    \]
  \item If $k_{\cX}=0$, then 
    \[
    I_{\cX}(t)= \iota_\star \(F(t)z+G(t)+O(z^{-1})\)
    \]
    for some functions $F\colon \C \to \Lambda$, $G \colon \C \to \HorbXL$, and 
    \[
    \iota_\star \, \bJ_{\cX} \(\tau(t)\)=\frac{I_{\cX}(t)}{F(t)}
    \quad
    \text{where}
    \quad
    \tau (t)=\frac{G(t)}{F(t)}.
    \]
  \end{enumerate}
\end{cor}

\begin{proof}
  The assertions $I_{\cX}(t)= \iota_\star (\cdots)$ follow by
  expanding $I_{\cX}(t)$ as a Laurent series in $z^{-1}$ and applying
  Lemma~\ref{6lem:1}.  The rest follows by combining
  Theorem~\ref{1thm:smallJ} with the ``Quantum Lefschetz'' theorem
  \cite{CCIT:Lefschetz}*{Corollary~5.1}.
\end{proof}

Corollary~\ref{1cor:ci_specific} follows immediately from
Corollary~\ref{6cor:hyp}, by computing the functions $\tau(t)$ in (1)
and $G(t)$ in (2) using Lemma~\ref{6lem:1}.

\begin{proof}[Proof of Proposition~\ref{1pro:terminal_singularities}]
  We recall the Reid--Tai criterion for terminal singularities \cite{YPG}.
  Fix a positive integer $r$ and a set of integer weights $a_1,\dots,
  a_n$ and consider the space
\[
\frac1{r}(a_1,\dots,a_n):=\C^n /\mu_r \quad \text{where}\; \mu_r \;
\text{acts with weights}\; a_1, \dots, a_n.
\]
We say that the set of weights is \emph{well-formed} if $\hcf (r, a_1,
\dots, \widehat{a_i}, \dots, a_n)=1$ for all $i$, that is if the
action of $\mu_r$ is faithful and there are no quasi-reflections. This
means that the orbifold is ``nonsingular'' in codimensions $0$ and
$1$.  The Reid--Tai criterion states that $X$ is well-formed with
terminal singularities if and only if
\begin{equation}
  \label{6eq:reid_tai}
  \sum_{i=1}^n \Big\langle \frac{ka_i}{r}\Big \rangle > 1 \quad
  \text{for  $k=1,2,\dots,r-1.$}  
\end{equation}
Terminal singularities are defined in \cite{YPG}; for the purpose of
this proof, the reader can take the Reid--Tai criterion as a
definition.

We now proceed to the proof of the Proposition. Let us assume that
$\cX=\cX_{d_0, \dots, d_m}\subset \p^\mathbf{w}$ is quasismooth and
well-formed with terminal singularities. Choose a non-zero $f\in F$.
Assuming that
\[
c=\# \{i\mid fw_i \in \z\}-\#\{j\mid fd_j \in \z \}\geq 0,
\] 
we want to show that 
\begin{equation}
  \label{6eq:condition_box}
  \sum_{i=0}^n \langle fw_i\rangle > 1+ \sum_{j=0}^m \langle d_jf\rangle.
\end{equation}
As in the proof of Lemma~\ref{6lem:1}, we can reorder the $d_j$ and
the $w_i$ so that:
\begin{enumerate}
\item for $j\leq l$, $fd_j \equiv fw_j \mod \z$ and none of these
  numbers is an integer.
\item $fd_j\in \z$ for $l<j$ and $fw_i\in \z$ for $l<i\leq m+c$. 
\end{enumerate}
The singularities of $\cX$ along $\p(V^f)$ are locally of the form:
\begin{equation}
  \label{6eq:singularity}
\frac1{r}\Big(\underbrace{0,0,\ldots,0}_c,w_{m+c+1},\dots,w_n\Big)
\end{equation}
Equation~\eqref{6eq:condition_box} is equivalent to
\[
\sum_{i=m+c+1}^n \langle fw_i\rangle >1
\]
and it holds by the Reid--Tai criterion for the singularity
\eqref{6eq:singularity}. The above argument can be read in reverse to
show the converse: if the condition of
Proposition~\ref{1pro:terminal_singularities} holds, then $\cX$ has
terminal singularities.
\end{proof}

\end{document}